\documentclass[onefignum,onetabnum]{siamart171218}

\usepackage{multirow}
\usepackage{amsfonts,amsmath,amssymb}
\usepackage{graphicx}
\usepackage{units}
\usepackage{graphicx}
\usepackage{url}
\usepackage{textcomp}
\usepackage{array}
\ifpdf
\DeclareGraphicsExtensions{.eps,.pdf,.png,.jpg}
\else
\DeclareGraphicsExtensions{.eps}
\fi

\usepackage[right=3cm,left=3cm]{geometry}

\newsiamremark{remark}{Remark}
\newsiamremark{hypothesis}{Hypothesis}
\crefname{hypothesis}{Hypothesis}{Hypotheses}
\newsiamthm{claim}{Claim}
\newsiamthm{preposition}{Preposition}

\headers{Node generation for meshless discretizations}{J. Slak and G. Kosec}

\title{On generation of node distributions for meshless PDE discretizations
  \thanks{Submitted to the editor on December 7th, 2018. Accepted August 19th, 2019.
    To appear in SIAM Journal on Scientific Computing.
    \funding{This work was supported by FWO grant G018916N, the ARRS research core funding no.~P2-0095 and Young Researcher program PR-08346.}}}

\author{Jure Slak\footnotemark[3] \thanks{Faculty of Mathematics and Physics, University of Ljubljana, Jadranska 19, 1000 Ljubljana, Slovenia
    (\email{jure.slak@ijs.si}, \url{http://e6.ijs.si/\string~jslak/}).}
  \and Gregor Kosec\thanks{``Jožef Stefan'' Institute, Department E6, Parallel and Distributed Systems Laboratory, Jamova cesta 39, 1000 Ljubljana, Slovenia
    (\email{gregor.kosec@ijs.si}, \url{http://e6.ijs.si/\string~gkosec/}).}
}

\usepackage{amsopn}

\ifpdf
\hypersetup{
  pdftitle={On generation of node distributions for meshless PDE discretizations},
  pdfauthor={J. Slak and G. Kosec}
}
\fi

\usepackage{algpseudocode}  
\usepackage{algorithm}      
\algnewcommand\algorithmicto{\textbf{to}}
\algnewcommand\algorithmicin{\textbf{in}}
\algnewcommand\algorithmicforeach{\textbf{for each}}
\algrenewtext{For}[3]{\algorithmicfor\ #1 $\gets$ #2\ \algorithmicto\ #3\ \algorithmicdo}
\algdef{S}[FOR]{ForEach}[2]{\algorithmicforeach\ #1\ \algorithmicin\ #2\ \algorithmicdo}

\newcommand{\R}{\mathbb{R}}
\newcommand{\N}{\mathbb{N}}

\renewcommand{\b}{\boldsymbol}
\newcommand{\eps}{\varepsilon}
\renewcommand{\phi}{\varphi}
\newcommand{\dpar}[2]{\frac{\partial #1}{\partial #2}}
\newcommand{\dr}{\ensuremath{h}}
\newcommand{\xm}{\ensuremath{x_\text{min}}}
\newcommand{\xM}{\ensuremath{x_\text{max}}}
\newcommand{\ym}{\ensuremath{y_\text{min}}}
\newcommand{\im}{\ensuremath{i_\text{min}}}
\newcommand{\yM}{\ensuremath{y_\text{max}}}
\newcommand{\X}{\mathcal{X}}
\DeclareMathOperator{\bb}{bb}
\DeclareMathOperator{\obb}{obb}

\begin{document}

\maketitle

\begin{abstract}
  In this paper we present an algorithm that is able to generate locally regular
  node layouts with spatially variable nodal density for interiors of arbitrary domains
  in two, three and higher dimensions. It is demonstrated that the generated node distributions
  are suitable to use in the RBF-FD method, which is demonstrated by solving
  thermo-fluid problem in 2D and 3D.
  Additionally, local minimal spacing guarantees
  are proven for both uniform and variable nodal densities.
  The presented algorithm has time complexity $O(N)$ to generate
  $N$ nodes with constant nodal spacing and $O(N \log N)$ to generate variably spaced nodes.
  Comparison with existing algorithms is performed
  in terms of node quality, time complexity, execution time and PDE solution accuracy.
\end{abstract}



\section{Introduction}
In recent years, a number of meshless approaches have been developed to numerically solve
partial differential equations (PDEs) with the desire to
circumvent the problem of polygonization encountered in the
classical mesh-based numerical methods. The major advantage of meshless methods is the ability to solve
PDEs on a set of scattered nodes, i.e.\ without a mesh. This
advantage was advertised even to the point that arbitrary nodes could be used
(see~\cite[p.~14]{liu2002mesh} and \cite{Reuther2012}), making node generation
seemingly trivial. Nevertheless, it soon turned out that such simplification leads to unstable
results.

Although node placing is considered much easier than mesh generation,
certain care still needs to be taken when generating node sets for meshless methods.
Many methods require sufficiently regular nodes for adequate precision and stability.
Among others, this also holds for the
popular Radial Basis Function-generated Finite Differences method (RBF-FD)~\cite{fornberg2015primer}.
Despite the need for quality node distributions,
solving PDEs with strong
form meshless methods utilizing radial basis functions (RBFs) has become
increasingly popular~\cite{Fornberg_Flyer_2015},
with recent uses in linear elasticity~\cite{slak2018refined},
contact problems~\cite{slak2018adaptive}, geosciences~\cite{fornberg2015primer},
fluid mechanics~\cite{kosec2018local}, dynamic thermal rating of
power lines~\cite{kosec2018afm} and even in the financial sector~\cite{golbabai2017new}.

Since one of the key advantages of mesh-free methods is the ability to use highly
spatially variable node distributions, which can adapt to irregular geometries and allow
for refinement in critical areas, many specialized algorithms for generations of such node
layouts have been developed. Most of them can generally be categorized into either
mesh-based, iterative, advancing front or sphere-packing algorithms.

The most basic way to generate such node sets is to employ existing tools and algorithms
for mesh generation, use the generated nodes and simply discard the connectivity
relations. Such approach was reasoned by Liu~\cite[p.~14]{liu2002mesh} as:
``There are very few dedicated node generators available commercially; thus, we
have to use preprocessors that have been developed for FEM.''.
This is problematic for two reasons: it is computationally wasteful, and some
authors have reported that such node layouts yielded unstable operator
approximations~\cite{shankar2018robust}, making them unable to obtain a solution.
Besides above shortcomings, such approach is also conceptually flawed, since the
purpose of mesh-free methods is to remove meshing from the solution procedure altogether.
To this end, other approaches were researched, often inspired by the algorithms
for mesh generation.

A common iterative approach is to position nodes by simulating free charged
particles, obtaining so-called minimal energy nodes~\cite{hardin2004discretizing}.
Other iterative methods include bubble simulation~\cite{liu2010node},
Voronoi relaxation~\cite{Balzer_2009} or a combination of both~\cite{choi1999node}.
Iterative methods are computationally expensive and require an initial distribution.
Additionally, the user is often required to consider trade-offs
between the number of iterations and node quality.
Despite their expensive nature, the produced distributions are often of high quality,
which makes iterative methods useful for improving node distributions generated by
other algorithms~\cite{fornberg2015fast}.

The next category consists of
advancing front methods, which usually begin at the boundary and advance towards
the domain interior, filling it in the process.
Löhner and Oñate~\cite{lohner2004general} present a general advancing front technique
that can be used for filling space with arbitrary objects. These methods, especially if
generating a mesh, are often restricted to two dimensions~\cite{persson2004simple}.
Another example of a two-dimensional
advancing front approach is inspired by dropping variable-sized grains into
a bucket~\cite{fornberg2015fast}, which yields quality variable density node distributions
and is computationally efficient in practice~\cite{slak2018fast}.

The last category are the circle or sphere packing methods~\cite{li2000point},
which generate high quality node distributions, but are often computationally expensive.
With inspiration from the graphics community,
Poisson disk sampling~\cite{cook1986stochastic} has become of interest. It
can be used to efficiently generate nodes in arbitrary dimensions~\cite{bridson2007fast},
and has just recently been used as a node generation
algorithm~\cite{shankar2018robust} providing nodal distributions of sufficient
quality for the RBF-FD method.

To the best of our knowledge, algorithms presented in~\cite{fornberg2015fast, shankar2018robust}
are currently among the best available. However, they have some shortcomings,
namely~\cite{fornberg2015fast} only works in two dimensions and~\cite{shankar2018robust} does not
support variable nodal spacing.
In this paper, we present an algorithm that overcomes these shortcomings.
The presented algorithm works in two, three and higher dimensions and supports
variable density distributions. It also has minimal spacing guarantees and
is provably computationally efficient. The main shortcoming of the presented
algorithm is that it requires discretized boundary as an input, which will be
addressed in future work. For algorithms that can fill domains with varying density,
conformal mappings can be used
to generate nodes on curved surfaces by appropriately modifying the node
density~\cite{fornberg2015fast}. The paper by Shankar et al.~\cite{shankar2018robust}
also includes an algorithm for generation of an appropriate boundary discretization,
based on RBF geometric model and super-sampling. This paper does not deal with the task of
generating a boundary discretization and focuses on discretizations of domain interiors,
assuming that the boundary discretization already exists when required. The extension
of the algorithm to curved surfaces will be addressed in future work.

The rest of the paper is organized as follows: in~\cref{sec:requirements} the
requirements for node generation algorithms are discussed, in~\cref{sec:known}
recently introduced algorithms for generating nodal distributions, suitable for strong
form meshless methods, are presented, in~\cref{sec:pds} a new algorithm is presented, in~\cref{sec:comp}
the algorithms are compared, and some numerical examples are presented in~\cref{sec:num}.

\section{Node placing algorithm requirements}
\label{sec:requirements}

In this section we examine a list of properties that an ideal node-positioning algorithm should
possess and discuss the rationale behind each property. The properties are
loosely ordered by decreasing importance.

\begin{enumerate}
  \item \emph{Local regularity.}
  Nodal distributions produced by the algorithm should be locally regular throughout the domain,
  i.e.\ the distances between nodes should be approximately equal.
  This definition of local regularity is somewhat
  soft and imprecise. The requirement stems from the fact
  that local strong form meshless methods are often sensitive to nodal positions and
  large discrepancies in distances to the nearest neighbors or other irregularities
  can cause ill-conditioned approximations, making the distribution inappropriate
  for solving PDEs. Thus, this point should be read in practice as follows:
  ``The distributions produced by the algorithm should yield quality PDE solutions when using
  local strong form methods, if reasonable spacing function $h$ was given.''.
  \item \label{itm:min-dx} \emph{Minimal spacing guarantees.}
  Computational nodes that are positioned too closely
  can severely impact the stability of some meshless methods~\cite{liu2002mesh}. Thus,
  provable minimal spacing guarantees are desirable. For constant spacing $h$,
  the condition
  \begin{equation}
  \|p - q\| \geq h
  \end{equation}
  is required for any two different points $p$ and $q$.
  For variable nodal spacing, the algorithm should guarantee a local
  lower bound for internodal spacing.

  \item \label{itm:h-var} \label{itm:first} \emph{Spatially variable densities.}
  Many node distribution algorithms
  rely on a constant discretization step $h$, as do some efficient implementations~\cite{bridson2007fast}.
  Spatially variable
  nodal distributions are often required when dealing with irregular domains or adaptivity~\cite{slak2018adaptive}.
  The algorithm should be able to generate distributions with spatially variable nodal spacing,
  which can be assumed to be given as a function $h\colon\R^d\to(0, \infty)$. The changes in
  variability should be gradual and smooth in order to satisfy the requirement of local regularity.
  The algorithm should work without any continuity assumptions for reasonable $h$
  (see remarks in~\cref{sec:pds-remarks}) and should see a constant step $h$ as a special
  case of variable step $h(p)$, not the other way around.

  \item  \label{itm:time} \emph{Computational efficiency and scalability.}
  Time complexity of the algorithm should ideally be linear in number of generated nodes.
  Quasilinear time complexity (e.g.\ $O(N \log N)$) is acceptable, while time complexity
  that is $\Omega(N^\alpha)$, for $\alpha > 1$, is undesirable.
  The algorithm should also be computationally efficient in practice, making it feasible to use
  as a node generation algorithm in an adaptive setting.

  \item \emph{Compatibility with boundary discretization.}
  Assume that a boundary discretization $\X_b$ of $\partial \Omega$ conforming to the spacing function $h$,
  already exists. More precisely, we are given a set $\X_b$ of points such that for any two neighboring
  points $p$ and $q$, the norm $\|p-q\|$ is approximately equal to $h(p)$ or $h(q)$. The generated
  discretization of the whole domain $\Omega$ should seamlessly join with the boundary discretization.
  This helps to prevent problems often encountered when enforcing boundary conditions
  (see~\cite[sec.~3.5]{shankar2018robust} and references therein).

  \item \label{itm:irreg} \emph{Compatibility with irregular domains.} The algorithm should inherently
  work with any irregular domain $\Omega$, given its characteristic (i.e.\ ``is element of'')
  function
  \begin{align}
  \chi_\Omega&\colon\R^d\to \{0, 1\}, \nonumber \\
  \chi_\Omega(p)&= \begin{cases} 1, & p \in \Omega, \\
  0, & p \notin \Omega. \end{cases}
  \end{align}
  Any algorithms that fill axis- or otherwise oriented bounding boxes,
  or produce nodes in a certain non-constant space outside $\Omega$
  are seen as impaired in this aspect. Desirably,
  as the volume of $\Omega$ decreases, no matter what the shape of $\Omega$ is,
  so should the number of operations required by the algorithm.

  \item \label{itm:dim-ind} \emph{Dimension independence.}
  The algorithm should ideally be formulated for a general (low) dimension $d$ without any special cases.
  One-, two- and three- dimensional versions of the algorithm should also allow a single
  general implementation.

  \item \label{itm:dir-ind} \emph{Direction independence.}
  The produced distributions and running time of the algorithm
  should be independent of the orientation of the domain $\Omega$ or
  the coordinate system used.

  \item \label{itm:free-param} \emph{No free parameters.} The algorithm should aim to minimize the
  number of free or tuning parameters and work well for all domains and density functions, without
  any user intervention. The aim is to require algorithms to be robust and work ``out of the box''.
  Any free parameters should be explored and well understood, default values should be recommended,
  and varying the parameters should not drastically change the algorithm's behavior.

  \item \label{itm:impl} \label{itm:last} \emph{Simplicity.}
  Algorithms with simpler formulations and implementations are preferred.
\end{enumerate}

\section{State of the art algorithms}
\label{sec:known}

To the best of our knowledge, recently published algorithms by Fornberg and
Flyer~\cite{fornberg2015fast} and Shankar et
al.~\cite{shankar2018robust} best satisfy the requirements described in~\cref{sec:requirements} and
are hence used as a base for further development. Both algorithms are first briefly described in the
following sections.

\subsection{Algorithm by Fornberg and Flyer}
\label{sec:ff}
The node positioning algorithm by Fornberg and Flyer~\cite{fornberg2015fast} was
published in 2015 in a paper titled ``Fast generation of 2-D node distributions for
mesh-free PDE discretizations''. The algorithm in its base form
constructs discretizations for two dimensional
axis-aligned rectangles and is presented as~\cref{alg:ff}. In the following
discussion, the first letters of the authors' surnames (FF) will be used to refer to the algorithm.

\begin{algorithm}
    \small
  \caption{Node positioning algorithm by Fornberg and Flyer.}
  \label{alg:ff}
  \vspace{5pt}
  \textbf{Input:} Box $[\xm, \xM] \times [\ym, \yM]$, a function $\dr\colon[\xm, \xM] \times [\ym, \yM] \to (0, \infty)$ and $n \in \N$. \\
  \textbf{Output:} An array of points in $[\xm, \xM] \times [\ym, \yM]$ distributed according to $\dr$.
  \begin{algorithmic}[1]
    \Function{FF}{$\xm$, $\xM$, $\ym$, $\yM$, $\dr$, $n$}
    \State $pts \gets$ an empty array of points \Comment{This is the final array of points.}
    \State $candidates \gets$ points spaced according to $\dr$ from $\xm$ to $\xM$ at $y$ coordinate $\ym$
    \Comment{This variable represents potential point locations, candidates for actual points that will be in the final result.}
    \State $(\ym, \im) \gets \Call{findMinimum}{candidates}$ \Comment{Find minimal point with respect to $y$}
    \While{$\ym \leq \yM$} \hfill coordinate and return its value and index. \label{ln:ff:while}
    \State $p \gets candidates[\im]$ \label{ln:ff:c-start}
    \State append $p$ to $pts$ \label{ln:ff:c-end}
    \State remove points closer than $\dr(p)$ from $candidates$ \label{ln:ff:m-start} \label{ln:ff:remove}
    \State find nearest remaining points in $candidates$ to the left and to the right of $p$ \label{ln:ff:find}
    \State add $n$ new points to $candidates$, equispaced on the circular sector with center $p$,
    spanning from the nearest left to the nearest right point
    \State $(\ym, \im) \gets \Call{findMinimum}{candidates}$  \label{ln:ff:m-end}
    \EndWhile
    \State \Return $pts$
    \EndFunction
  \end{algorithmic}
\end{algorithm}

Initially, the lower side of the rectangle is filled with nodes, spaced according to the given
spacing function $\dr$. The algorithm works as an advancing front algorithm,
starting from $y = y_\text{min}$ and advancing towards $y = y_\text{max}$.
In each iteration the lowest (min$(y)$) candidate $p$ from the current list of
potential node locations is found, removed from potential candidates
and added to the final list. New candidates are spaced accordingly away
from $p$ and are inserted into the list of potential node locations.
The iteration continues until $y=\yM$ limit is reached.

For irregular domains $\Omega$ the authors recommend to run the above algorithm for the
bounding box of $\Omega$, denoted $\bb(\Omega)$, and later discard the nodes outside $\Omega$.
If present, the boundary discretization is superimposed
onto the discretization generated by~\cref{alg:ff}.
Additionally, internal nodes whose closest boundary node $p$ is less than $h(p)/2$ away
are discarded as well.

A few local ``repel algorithm'' iterations are recommended in the vicinity of
the boundary to improve the quality. This part will be omitted from consideration,
as it is an iterative improvement scheme that can be performed equivalently on any
node distribution generated by any other algorithm~\cite{kosec2018local}.
The behavior of FF near the boundaries is thus excluded from analysis, as
it is designed to work with the ``repel algorithm''.

\subsubsection{Time complexity analysis}
\label{sec:ff-time}
The complexity of the~\cref{alg:ff} is not given by the authors
and is hence derived in this section.
Denote the number of generated nodes with $N$ and the size of array $candidates$ at the start of
iteration $i$ with $s_i$.  Everything up to while loop on line~\ref{ln:ff:while}
is negligible compared to the main loop and takes $O(1)$ time for creation of lists
and $O(s_0)$ for candidate generation and minimum extraction. In the main loop,
lines~\ref{ln:ff:c-start}--\ref{ln:ff:c-end} consume (amortized) constant time
and lines~\ref{ln:ff:m-start}--\ref{ln:ff:m-end} take time proportional to the size
of $candidates$ array, i.e.\ $O(s_i)$ time. Total time complexity is therefore
\begin{equation}
T_{\text{FFbox}} = O(1) + O(s_0) + \sum_{i=1}^N (O(1) + O(s_i)) = O(N \max_{1 \leq i \leq N} s_i) := O(NS),
\end{equation}
where $S$ is defined as $S = \max_{1 \leq i \leq N} s_i$.

Precisely analyzing $s_i$ and $S$ is difficult for general function $\dr$. However,
for a fixed square box and constant spacing $h$ it holds that $N = \Theta(\frac{1}{\dr^2})$
and $s_i = \Theta(n \frac{1}{h}) = \Theta(n \sqrt{N})$.
The time complexity in this case is hence $O(nN\sqrt{N})$.


For irregular domains $\Omega$ additional work is required. If $N$ denotes the
final number of nodes, the algorithm will generate
approximately $\frac{|\bb{\Omega}|}{|\Omega|}N$ nodes in case of constant $\dr$.
Superimposing the boundary discretization with $N_b$ nodes and testing all generated nodes for proximity
takes $O(N_b\log N_b + \frac{|\bb{\Omega}|}{|\Omega|}N \log N_b)$ time, for building
and querying the $k$-d tree of boundary nodes. These terms are dominated by the node
generation in the interior and the time complexity of the algorithm by Fornberg and Flyer
for generating node distributions for irregular domains for constant spacing $\dr$ is
\begin{equation}
\label{eq:ff-time}
T_{\text{FF}} = O\left(n\left(\frac{|\bb{\Omega}|}{|\Omega|}N\right)^{1.5}\right).
\end{equation}
For variable spacing, the overhead of generated nodes due to the irregularity of
$\Omega$ and the advancing front size have to be evaluated using integrals, making the
time complexity expression somewhat more complicated and less illustrative.
\subsubsection{Implementation notes}
Authors provided a Matlab implementation of~\cref{alg:pds} in the
Appendix ot their article~\cite{fornberg2015fast}. This implementation has been translated to C++
using the Eigen matrix library~\cite{eigenweb} and the \texttt{nanoflann} library for $k$-d
trees, provided by Blanco and Rai~\cite{blanco2014nanoflann}.
The translation is mostly faithful to the original
with a few inefficiencies removed. The C++ implementation is approximately 6 times faster than the original Matlab implementation (both tested on the same computer).

\subsection{Algorithm by Shankar, Kirby and Fogelson}
\label{sec:skf}
In 2018, Shankar, Kirby and Fogelson published a node generation algorithm
in a paper titled ``Robust node generation for meshfree discretizations on irregular domains
and surfaces''~\cite{shankar2018robust}. Their node generation algorithm is designed to work on
surfaces and in 3-D, however it does not support variable nodal spacing. We will focus our attention
on the part that generates discretizations of the domain interior, when boundary discretization has
already been constructed. The main part of the node generation for the interior is based on
Poisson disk sampling of the oriented bounding box $\obb(\Omega)$ of domain $\Omega$,
described in a paper by Bridson~\cite{bridson2007fast}. The relevant part of the node generation
algorithm is presented as~\cref{alg:skf}. In the following discussion the first letters of the
authors' surnames (SKF) will be used to refer to the algorithm.

\begin{algorithm}
    \small
  \caption{Node positioning algorithm by Shankar, Kirby and Fogelson.}
  \label{alg:skf}
  \vspace{5pt}
  \textbf{Input:} Domain $\Omega$ and its dimension $d$. \\
  \textbf{Input:} A nodal spacing step $\dr > 0$. \\
  \textbf{Input:} A list of boundary points $\X$ of size $N_b$, moved $h$ towards domain interior. \\
  \textbf{Output:} A list of points in $\Omega$ distributed according to spacing function $\dr$.
  \begin{algorithmic}[1]
    \Function{SKF}{$\Omega$, $h$, $\X$, $n$}
    \State $obb \gets \Call{OBB}{\X}$ \Comment{Generate an oriented bounding box of $\X$ using PCA.}
    \State $G \gets d$-dimensional grid of size $h/\sqrt{d}$ of $-1$.  \Comment{Maps points to their indices.}  \label{ln:skf-pds-start}
    \State $p \gets$ uniform random node inside $obb$
    \State $G[\Call{index}{p}] \gets 0$  \Comment{\textsc{index} returns $d$-d index of $p$, and its sequential index is 0.}
    \State $S \gets [\,p\,]$  \Comment{Resulting list of accepted samples.}
    \State $A \gets \{0\}$  \Comment{Set of active indices.}
    \While{$A \neq \emptyset$}
    \State $i \gets \Call{randomElement}{A}$ \Comment{Get a uniform random element of $A$.}
    \State $b \gets \texttt{false}$  \Comment{Indicates if any valid points were generated.}
    \For{$j$}{$1$}{$n$}
    \State $p \gets$ uniform random point in annulus with center $S[i]$ and radii $h$ and $2h$
    \If{\textbf{not} $\Call{outside}{p, obb}$ \textbf{and not} $\Call{tooClose}{p, h, G, S}$}
    \State \Call{Add}{$A$, \textsc{size}($S$)}  \Comment{Add sequential index of $p$ to active set $A$.}
    \State $G[\Call{index}{p}] \gets \Call{size}{S}$  \Comment{Mark grid cell taken by $p$ as occupied.}
    \State $\Call{Append}{S, p}$  \Comment{Append $p$ to the list of accepted samples.}
    \State $b \gets \texttt{true}$  \Comment{Flag that an accepted sample was generated.}
    \State \textbf{break for} \label{ln:skf-break}
    \EndIf
    \EndFor
    \If{$b = \texttt{false}$}  \Comment{Point $S[i]$ failed to generate any accepted samples.}
    \State \Call{remove}{$A$, $i$}  \Comment{Point with index $i$ is removed from the active set.}
    \EndIf
    \EndWhile  \label{ln:skf-pds-end}
    \State $T \gets \Call{kdtreeInit}{\X}$  \Comment{Initialize spatial search structure on points $\X$.}
    \State $S \gets \Call{Filter}{T, S}$ \Comment{Discard points outside the region, bounded by $\X$.}
    \State \Return $S$
    \EndFunction
  \end{algorithmic}
\end{algorithm}

The algorithm starts by taking points on the boundary and their corresponding outward unit normals
and shifting them towards the domain's interior by $h$.
An oriented bounding box (OBB) of the shifted boundary points is then constructed
using Principal Component Analysis (PCA)~\cite{jolliffe2011principal} as described by
Dimitrov et al.~\cite{dimitrov2006bounding}.
The main part of the algorithm, spanning lines~\ref{ln:skf-pds-start} to~\ref{ln:skf-pds-end},
is the Poisson disk sampling algorithm, which generates the internal discretization of the
oriented bounding box using a background grid $G$ as a spatial search structure.
Finally, points outside the domain, bounded by $\X$, are discarded. Here, a $k$-d tree is used to test all
candidates for inclusion by testing against the outward normal of their closest
boundary point. The remaining points along with the original boundary discretization
constitute the final discretization. As an inward-shifted array of points was used
to construct the internal discretization, spacing of at least $h$ is guaranteed.

\subsubsection{Time complexity analysis}
\label{sec:skf-time}
Authors themselves provide the time complexity analysis of the algorithm. Translated to our notation,
the running time of the interior fill algorithm is
\begin{equation}
\label{eq:skf-time}
T_{\text{SKF}} = O\left(n \frac{|\obb(\Omega)|}{|\Omega|} N\right).
\end{equation}
This represents the running time of the Poisson disk sampling. The PCA analysis and tree construction
are linear or log-linear in $N_b$ and are thus dominated by the Poisson disk sampling.

\subsubsection{Implementation notes}
\label{sec:skf-impl}
There is a small difference between the algorithm as described by Shankar et al.~\cite{shankar2018robust}
and Bridson~\cite{bridson2007fast}. The Bridson version generates up to $n$ candidates for
each point and stops as soon as one candidate is accepted. The version in the SKF algorithm
generates all $n$ points and adds all accepted candidates. \Cref{alg:skf} uses the original, Bridson
version and to obtain the SKF version, one needs to remove the break statement on line~\ref{ln:skf-break}.
Since the authors of SKF algorithm claim to use a faithful implementation of algorithm as presented by
Bridson and only list the algorithm for completeness, we decided to use
the Bridson version in our tests. All matrix and tensor operations were again implemented using
the Eigen matrix library and the $k$-d tree operations were implemented using \texttt{nanoflann}.

\section{New node placing algorithm}
\label{sec:pds}

From the discussion presented in~\cref{sec:known} it is clear that although state
of the art placing algorithms provide a solid spatial discretization methodology for
strong form meshless methods, there is still room for improvements, especially in the
generalization to higher dimensions, flexibility regarding variable nodal density,
and treatment of irregular domains. Improving upon the work of Fornberg and
Flyer~\cite{fornberg2015fast} and Shankar et al.~\cite{shankar2018robust},
we propose a new algorithm that overcomes some of limitations of FF and SKF algorithms.
We will refer to the proposed node placing algorithm as PNP in the rest of the paper.

The PNP algorithm is, similarly to SKF, based on Poisson disk sampling.
Poisson disk sampling has certain stochastic properties, such as the fact that it produces a
``blue noise distribution'' that is an excellent fit
for graphical applications like dithering~\cite{bridson2007fast,cook1986stochastic}.
In context of PDE solution procedure a slightly different
distribution is required that primarily follows appropriate spacing and regularity criteria.
Therefore, the PNP algorithm deviates from the original Poisson disk sampling
in order to effectively produce tightly packed regular distributions needed in solution of PDEs.

\begin{figure}[h]
  \centering
  \includegraphics[width=0.7\linewidth]{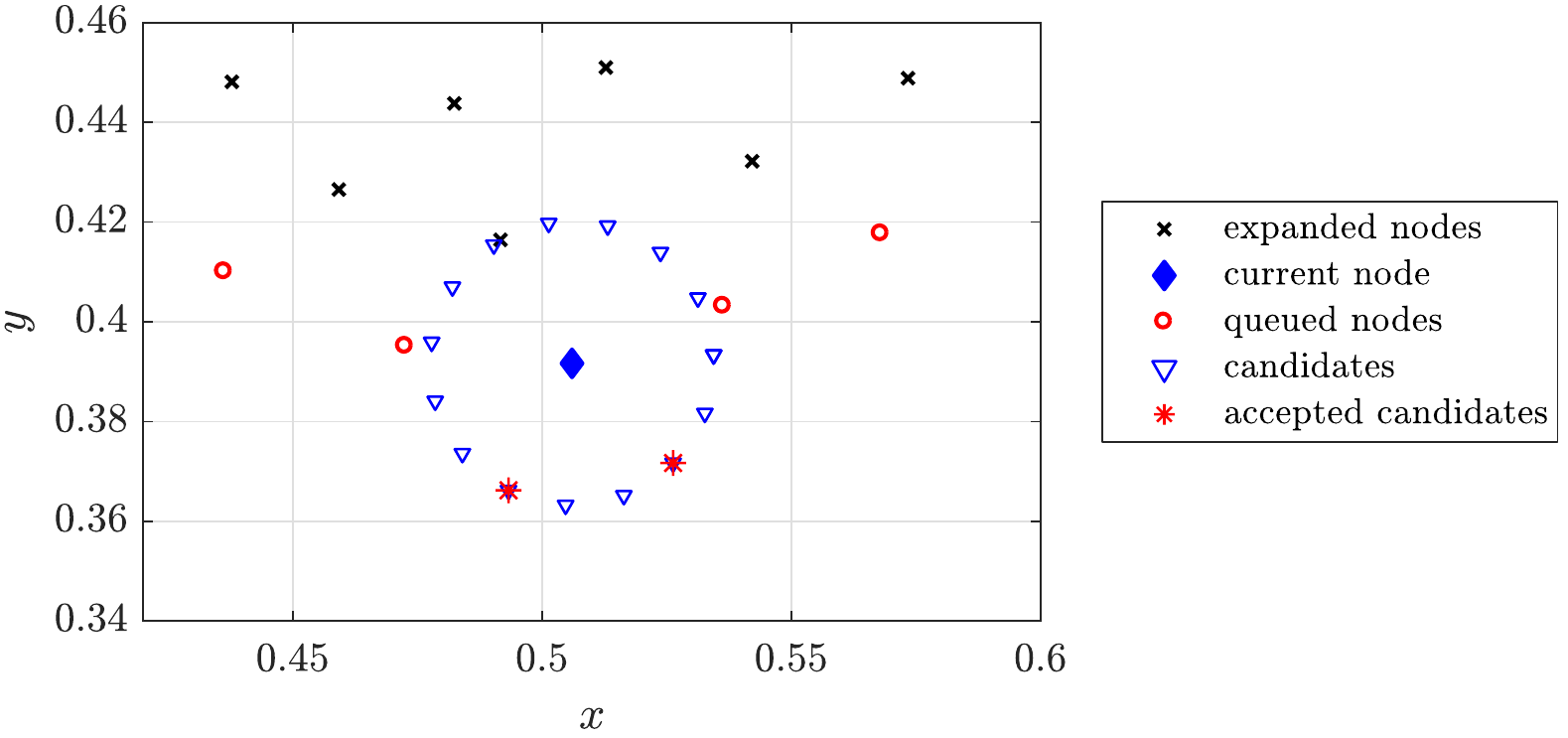}
  \caption{Generation and selection of new candidates in the proposed algorithm.}
  \label{fig:pds-new-candidates}
\end{figure}

The PNP algorithm begins either with a given non-empty set of ``seed nodes'' or
with an empty domain.
In the context of PDE discretizations, some nodes on the boundary are usually
already known and can be used as seed nodes, possibly along with additional nodes in the interior.
If algorithm starts with no nodes, it adds a seed node randomly within the domain.
Before the main iteration loop, seed nodes are put in a queue, waiting to be processed.
In each iteration $i$, a node $p_i$ is dequeued and \emph{expanded}, by generating a set of
candidates $C_i$, which are positioned on a sphere
with center $p_i$ and radius $r_i$, where $r_i$ is obtained from the function $h$, $r_i = h(p_i)$.
Candidates that lie outside of the domain or are too close to already existing nodes
are rejected and remaining candidates are enqueued for expansion.
Node $p_i$ is accepted as a domain node and will not be touched any more.
The iteration continues until the queue is empty. \Cref{fig:pds-new-candidates}
demonstrates a core operation of the algorithm, i.e.\ the expansion, with possible
selection of new candidates and the rejection process.

\begin{figure}[h]
  \centering
  \includegraphics[width=0.19\linewidth]{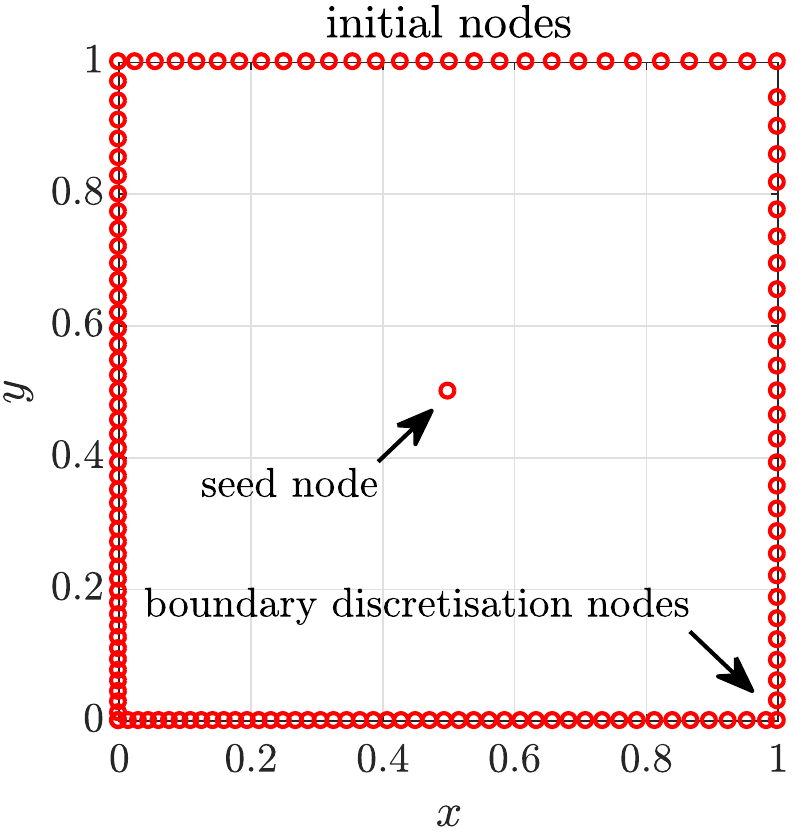}
  \includegraphics[width=0.19\linewidth]{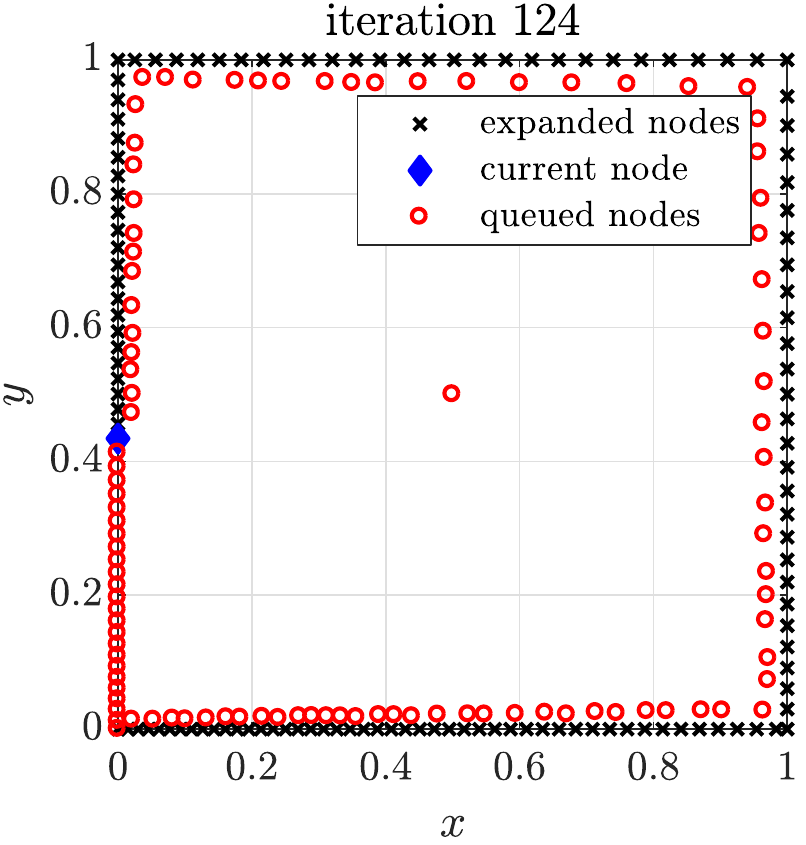}
  \includegraphics[width=0.19\linewidth]{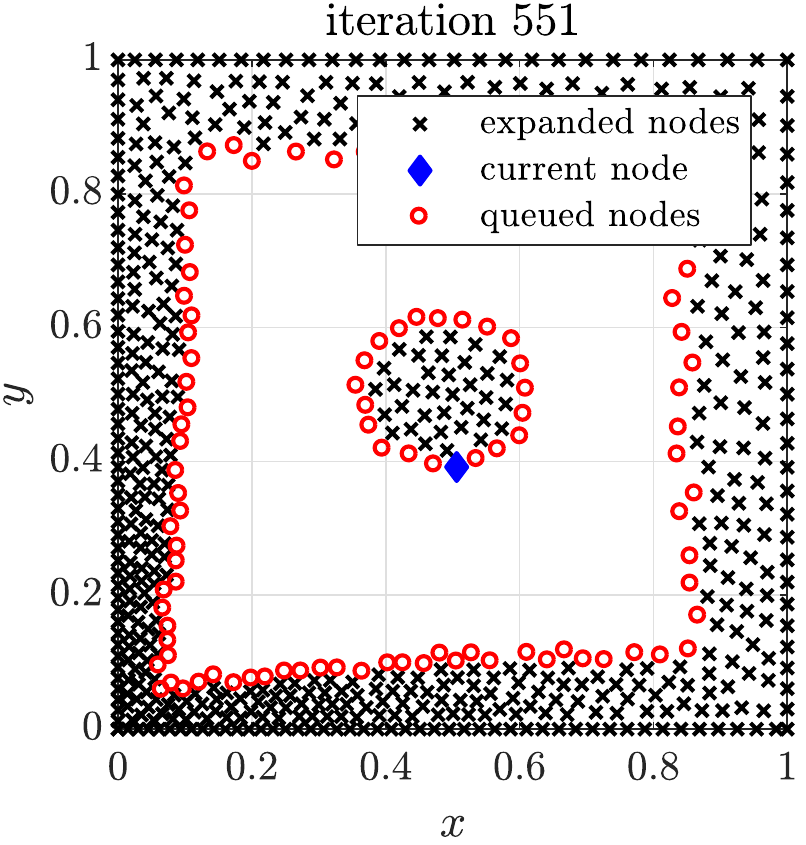}
  \includegraphics[width=0.19\linewidth]{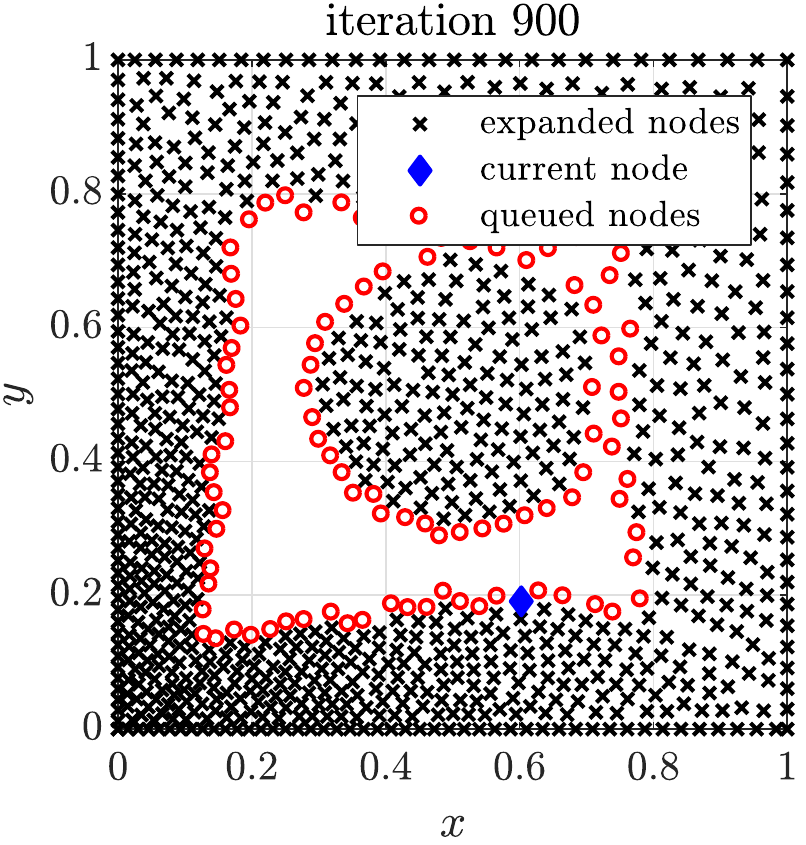}
  \includegraphics[width=0.19\linewidth]{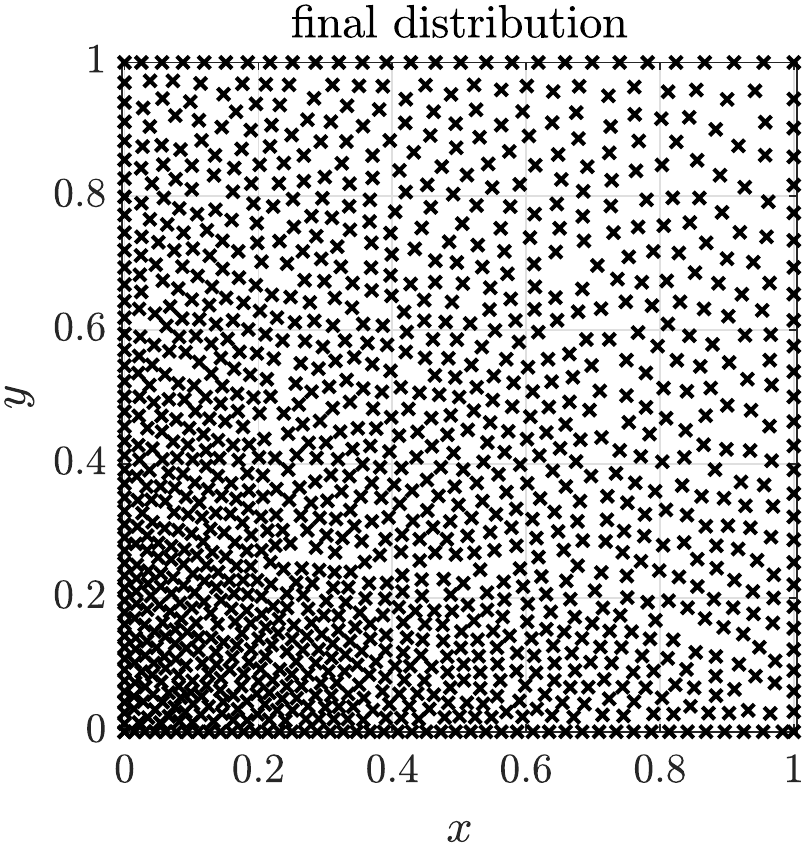}
  \caption{Run-time progress of the proposed algorithm (left to right).
    Unit square $[0, 1]^2$ was sampled with nodal spacing function $h(x, y) = 0.015\,(1 + x + y)$.}
  \label{fig:poisson-progress}
\end{figure}

\Cref{fig:poisson-progress} illustrates the execution of the algorithm.
The first panel shows an initial setup on a unit square. For demonstration purposes,
the nodes in the initial boundary discretization along with a single node in the
interior were chosen as the seed nodes.
The subsequent panels in~\cref{fig:poisson-progress} illustrate the progression
of the algorithm. The discretization grows from the initial nodes inwards towards
the empty interior, until no more acceptable candidates can be found due
to already existing nodes. The advancing front nature of the algorithm can be seen,
however the front itself is not a straight line as in FF.

An efficient implementation with an implicit queue contained in the array of final points and
the $k$-d tree spatial structure~\cite{moore1991intoductory}
is presented as~\cref{alg:pds}. In practice, the comparison
on line~\ref{ln:comapre-eps}
should be done with some tolerance, due to inexactness of the floating point arithmetic,
i.e.\ the line should in practice read $\|c_{i,j} - n_{i,j}\| \geq (1-\eps) r_i$ for e.g.\ $\eps = 10^{-10}$.

\begin{algorithm}
  \small
  \caption{Proposed node positioning algorithm.}
  \label{alg:pds}
  \vspace{5pt}
  \textbf{Input:} Domain $\Omega$ and its dimension $d$. \\
  \textbf{Input:} A nodal spacing function $\dr\colon\Omega \subset \R^d \to (0, \infty)$. \\
  \textbf{Input:} A list of starting points $\X$, this includes the possible boundary discretization and seed nodes. \\
  \textbf{Output:} A list of points in $\Omega$ distributed according to spacing function $\dr$.
  \begin{algorithmic}[1]
    \Function{PNP}{$\Omega$, $h$, $\X$}
    \State $T \gets \Call{kdtreeInit}{\X}$  \Comment{Initialize spatial search structure on points $\X$.}
    \State $i \gets 0$  \Comment{Current node index.}
    \While{$i < |\X|$}          \Comment{Until the queue is not empty.}
    \State $p_i \gets \X[i]$  \Comment{Dequeue current point.}
    \State $r_i \gets h(p_i)$  \Comment{Compute its nodal spacing.}
    \ForEach{$c_{i,j}$}{\Call{candidates}{$p_i, r_i$}} \Comment{Loop through candidates.} \label{ln:candidates}
    \If{$c_{i,j} \in \Omega$}    \Comment{Discard candidates outside the domain.}
    \State $n_{i,j} \gets \Call{kdtreeClosest}{T, c_{i,j}}$  \Comment{Find nearest node for
      proximity test.}
    \If{$\|c_{i,j} - n_{i,j}\| \geq r_i $}  \Comment{Test that $c_{i,j}$ is not too close to other nodes.}  \label{ln:comapre-eps}
    \State \Call{append}{$\X, c_{i,j}$}  \Comment{Enqueue $c_{i,j}$ as the last element of $\X$.}
    \label{ln:add-node}
    \State \Call{kdtreeInsert}{$T, c_{i,j}$}  \Comment{Insert $c_{i,j}$ into the spatial search structure.}
    \EndIf
    \EndIf
    \EndFor
    \State $i \gets i + 1$  \Comment{Move to the next non-expanded node.}
    \EndWhile
    \State \Return $\X$
    \EndFunction
  \end{algorithmic}
\end{algorithm}

The algorithm includes generation of new candidates in line~\ref{ln:candidates} that needs to be further defined. Three options are proposed and evaluated below:
\begin{enumerate}
  \item \textit{Random candidates:} The candidate set $C_i$ in each iteration consists of
  $n$ random points chosen on a $d$-dimensional sphere with center $p_i$ and radius $r_i$,
  reminiscing the original Poisson disk sampling.
  \item \textit{Fixed pattern candidates:} The candidate set $C_i$ in each iteration consists
  of a fixed discretization of a unit ball, translated to $p_i$ and scaled by $r_i$.  \label{item:fixed}
  The discretization of a unit ball in 2-D is obtained simply by
  $C_{\text{unit}}(k) = \{ (\cos\phi, \sin\phi); \ \phi \in \{0, \phi_0, 2\phi_0, \ldots, (k-1)\phi_0\}, \phi_0 = \frac{2\pi}{k} \}$.
  In $d$-dimensions, the discretization of a ball with radius $r$
  is obtained using $d$-dimensional spherical coordinates and recursively discretizing a $d-1$ dimensional ball.

  Using e.g.\ $k=6$ results in $14$ candidates in 3-D, and using $k=12$ results in $48$.
  In 3-D, the parameter $k$ represents the number of points lying on the great circle.

  \item \textit{Randomized pattern candidates:} The candidate set $C_i$ is obtained from
  the fixed set in point~\ref{item:fixed}, by applying a random rotation to all the points.
\end{enumerate}

\begin{figure}[h]
  \centering
  \includegraphics[width=0.24\linewidth]{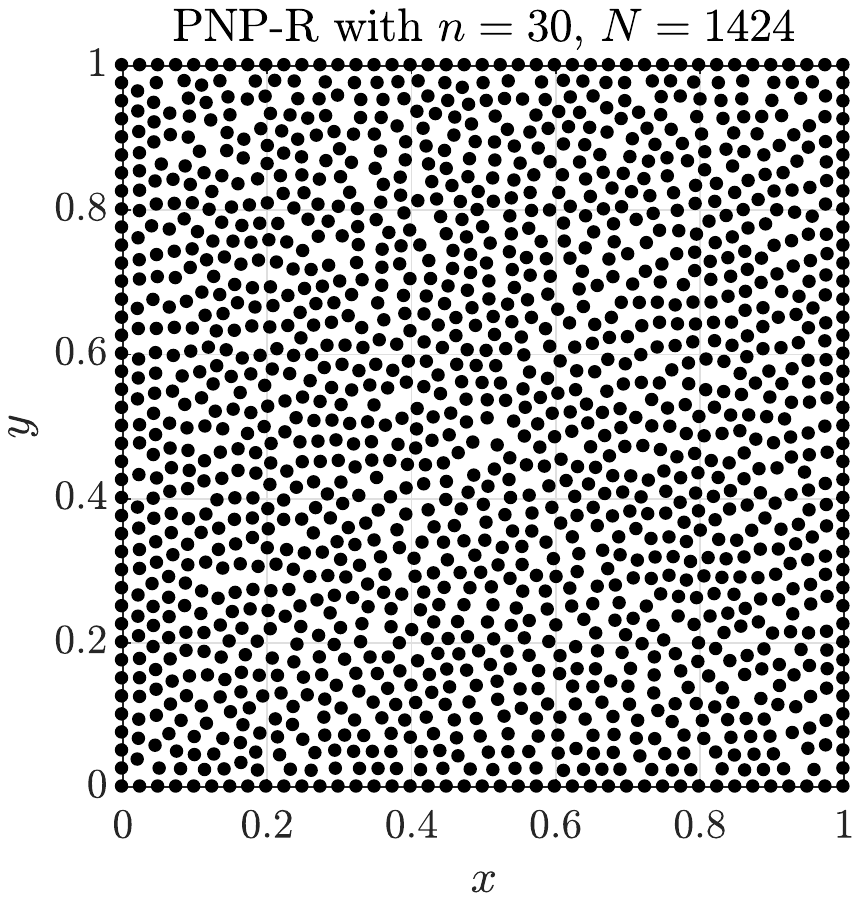}
  \includegraphics[width=0.24\linewidth]{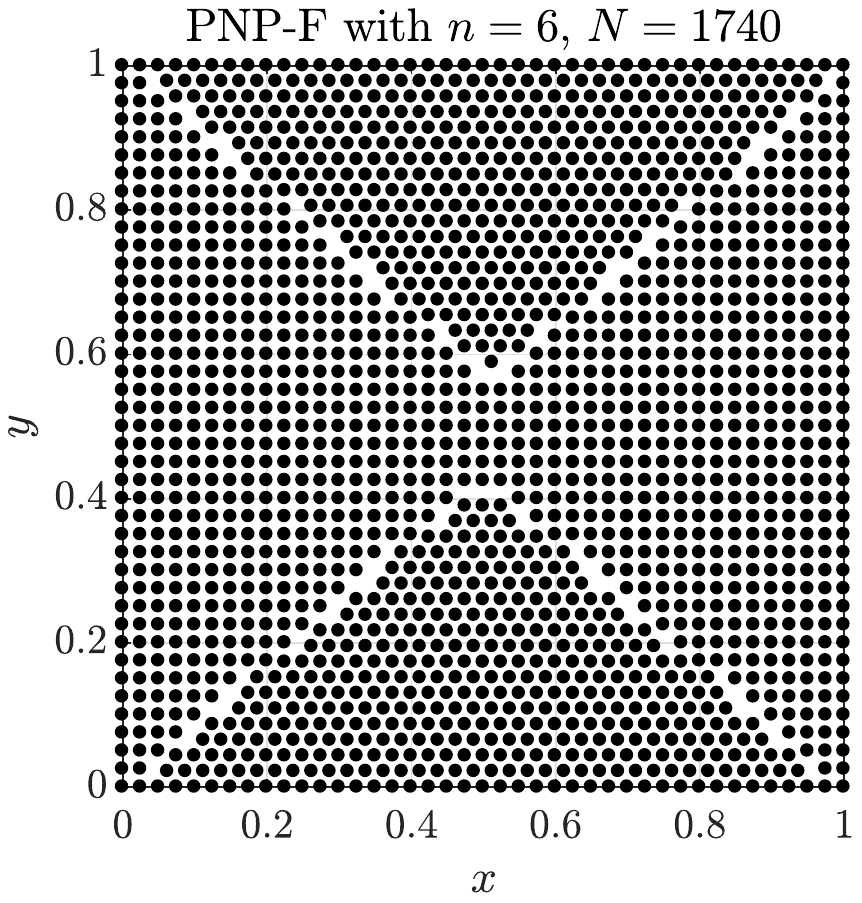}
  \includegraphics[width=0.24\linewidth]{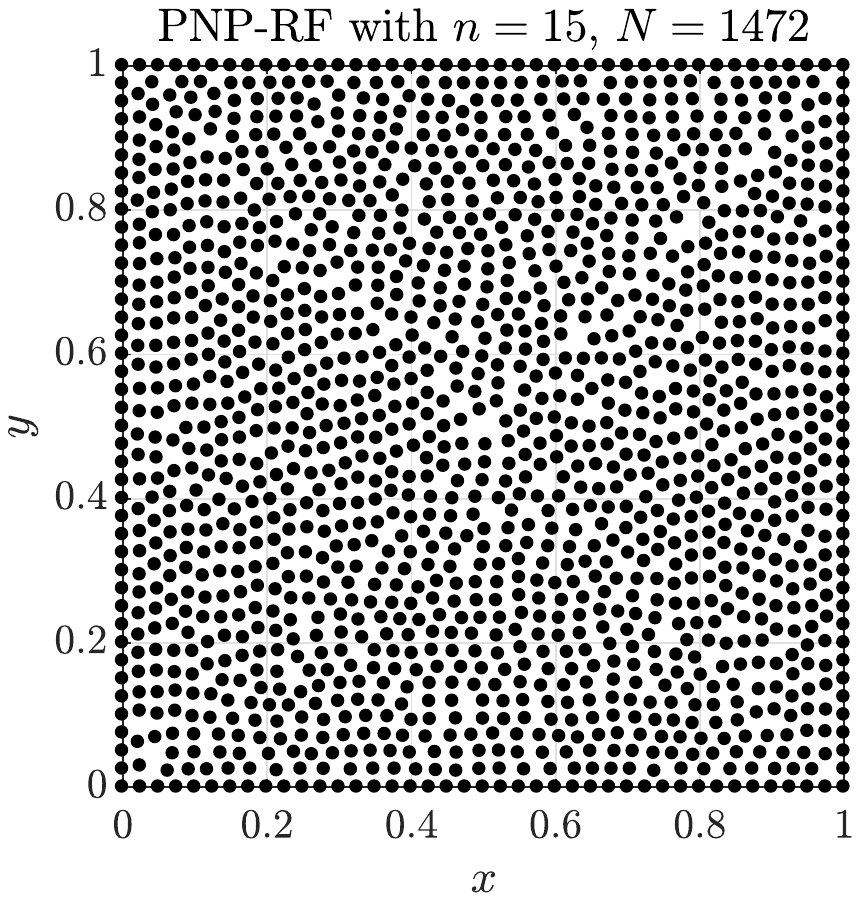}
  \caption{Comparison of different types of candidate generation when
    filling the unit square $[0, 1]^2$ with $h=0.025$.}
  \label{fig:pds-cand-choice}
\end{figure}

The three ways of candidate generation were used to produce node distributions on
a unit square, shown in~\cref{fig:pds-cand-choice}. Different
types of candidate generation are abbreviated as PNP-R, PNP-F and PNP-RF for random,
fixed pattern, and randomized fixed pattern variants, respectively.

The fixed pattern candidate generation algorithm stands out, as
the gaps where the advancing fronts of nodes joined are clearly visible.
Due to reproduction of space-efficient hex-packing it also has the most nodes.
Other two algorithms generate visually similar distributions, but a higher value of $n$
needed to be used for the randomized version to produce similar results.
We decided to use the randomized fixed pattern for candidate distribution,
as it produces results similar to the random version with lower time complexity.

The presented algorithm has a few differences compared to the original Poisson disk sampling.
First and foremost, the algorithm works with variable nodal spacing and is able to generate
distributions with spatially variable densities.
Each node is used only once to generate the new candidates.
Third, the candidates are generated uniformly on the sphere with random offsets, as opposed
to being generated at random on an annulus. This packs the candidates more tightly and
reduces the gaps. It also improves running time, as candidates better cover
the unoccupied space at the cost of removing the stochastic properties
of the sampling, which are not relevant to the PDE solutions.
Fourth, the candidates that are outside $\Omega$ are immediately discarded, only continuing
the generation of candidates actually inside $\Omega$, once again improving execution time.
More details about impact of above differences are investigated in~\cref{sec:comp}
and can be observed in~\cref{fig:sample-square} and~\cref{fig:time}.

PNP algorithm exhibits gaps between nodes where the advancing fronts meet in
\cref{fig:pds-cand-choice}, however the gaps are never large enough
that another node could be placed inside and are even emphasized visually is the marker size is
comparable to nodal spacing (see \cref{fig:sample-square}). The exact place where the
advancing fronts meet is dependent on the position of the seed nodes.
If the algorithm is run from a seed node in the domain interior instead of from the boundary
nodes, these types of front do not appear throughout the domain, but gaps form at the boundary
instead. With PDE discretization in mind it is less problematic to have them appear in the domain
interior, and this is how the algorithm is run for the rest of the paper.

Additionally, the algorithm can be easily modified to return the indices of the nodes where the
fronts meet. For each node $i$, we can check if any candidates generated from it are
accepted and added on line~\ref{ln:add-node}. If that is not the case, index $i$ can be
added to the list of \emph{terminal nodes}, which is returned after the algorithm finishes.
Regularization can then be performed on those or neighboring nodes, if necessary.

\subsection{Time complexity analysis}
\label{sec:pds-time-complexity}
Output sensitive time complexity is straightforward to analyze.
Let us denote the number of given starting points in $\X$ with $N_b = |\X|$.
Is is assumed that $N_b$ is significantly less than $N$, e.g.\
$N_b = O(N^\frac{d-1}{d})$ as is the case when $\X$ represents the
boundary discretization. The initial construction of the spatial index
costs $O(N_b \log N_b)$ and initialization of other variables costs $O(1)$.
The number of iterations of the main loop is equal to the number of generated
points, denoted by $N$. A total of $n$ candidates are generated in $i$-th iteration
and, in the worst case, two $k$-d tree operation on the tree with at most $i+N_b$ nodes are performed,
taking $O(\log (i+N_b))$ time for each candidate. All other operations are (amortized) constant.
Thus the total time complexity of the algorithm is equal to
\begin{equation}
T_{\text{PNP}} = O(N_b \log N_b) + O(1) + \sum_{i=1}^N \left[ n O(\log (i+N_b)) + O(1) \right] = O(nN\log N).
\end{equation}

The above analysis shows that the time complexity of the algorithm is dominated by the spatial
search structure used, which adds an undesired factor $\log N$.
If $h$ is assumed to be constant, the algorithm
could be sped up by using a uniform-grid based spatial search structure, similar to one used
in~\cref{alg:skf}. Using such a search structure requires a rectangular grid,
usually with spacing $h/\sqrt{d}$, such that there is at most one point per grid cell.
When constructed on the rectangle $\obb(\Omega)$, the time complexity of its allocation
and initialization is proportional to the number of cells, which leads to time complexity
\begin{equation}
O\left(\frac{|\obb\Omega|}{(h/\sqrt{d})^d}\right) =
O\left(\frac{|\obb \Omega|}{|\Omega|}N\right),
\end{equation}
using the fact that for constant $h$ the number of nodes is $N = \Theta(|\Omega|/h^d)$.

The subsequent insertions and queries in the grid are all $O(1)$, thus improving
the time complexity of the algorithm for constant $h$ to
\begin{equation}
T_{\text{PNP-grid}} = O\left(\frac{|\obb \Omega|}{|\Omega|}N + nN\right).
\end{equation}
Furthermore, the factor $\frac{|\obb \Omega|}{|\Omega|}$ can be eliminated by using a
hash map of cells instead of a grid; however, the practical benefit of that approach
shows only with very irregular domains.

Using the background grid for a spatial structure
is feasible even with moderately spatially variable $h$,
by allowing more than one point per cell. For even higher
variability, hierarchical grids could be used, but a $k$-d tree-like
search structure covers all cases. For a specific use case, $k$-d tree can
be replaced with any spatial search structure, as desired by the user,
obtaining time complexity
\begin{equation}
T_{\text{PNP-general}} = O(P(N) + Nn(Q(N)+I(N))),
\end{equation}
where $P$ is the precomputation/initialization time used by the data structure on $N$ nodes,
$Q(N)$ is the time spent on a radius query and $I(N)$ is the time spent for new element insertion.

\subsection{Implementation notes}
As in the previous two algorithms, all matrix and tensor operations were implemented using
the Eigen matrix library and the $k$-d tree operations were implemented using \texttt{nanoflann}.

\subsection{Remarks}
\label{sec:pds-remarks}
\Cref{alg:ff} and~\cref{alg:pds} do not necessarily terminate,
depending on the nodal spacing function used.
The integral
\begin{equation}
N(h) := \int_{\Omega} \frac{\textrm{d}\Omega}{h(p)^d},
\end{equation}
approximately measures the number of points required and can be infinite
even if function $h$ is smooth and positive on $\Omega$. Simply taking
a one dimensional example $\Omega = (0, 1)$ and $h(x) = \frac{0.1}{x}$
is enough to trick the algorithm into sampling indefinitely.
As a precaution to that and more practically, as a memory limit, the maximal
number of points $N_{\text{max}}$ can be specified by the user and the algorithm
can be terminated prematurely.

\section{Satisfaction of the requirements}
\label{sec:comp}

This section compares all three node placing algorithms, namely FF (\cref{alg:ff}),
SKF (\cref{alg:skf}) and PNP (\cref{alg:pds}).
The results of the comparison presented in this section are summarized at the end
in~\cref{tab:compare}.
The following subsections roughly follow the requirements postulated in~\cref{sec:requirements}.

\subsection{Local regularity}
\label{sec:node-quality}
The most important feature that an algorithm should possess is regularity of the distributions.
This property is initially tested visually, by observing plots of nodal distributions, which
is feasible only in 2-D. Among other things, local regularity states also that large discrepancies
in distances to nearest neighbors are not desired. This can be tested in arbitrary dimensions,
by observing distances to nearest neighbors, using various statistics and histogram plots to
determine their properties. Finally, accuracy and stability of
solutions of PDEs on generated node distributions can be compared to fully determine the quality
of distributions generated by the three algorithms.

\begin{figure}[h]
  \centering
  \includegraphics[width=0.24\linewidth]{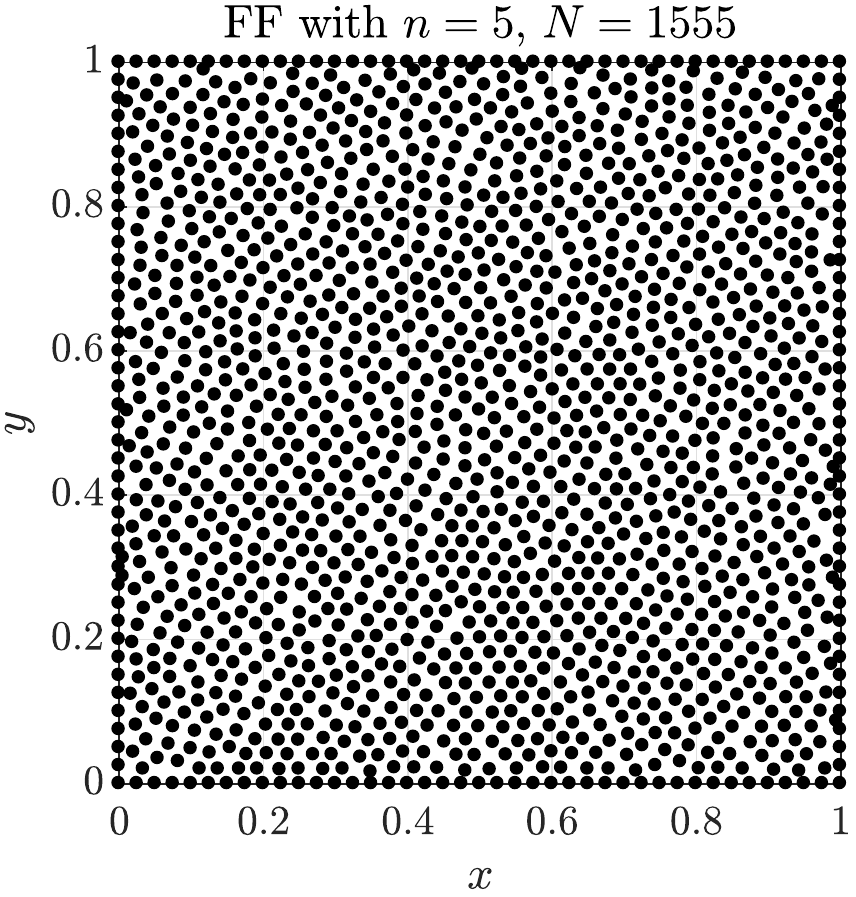}
  \includegraphics[width=0.24\linewidth]{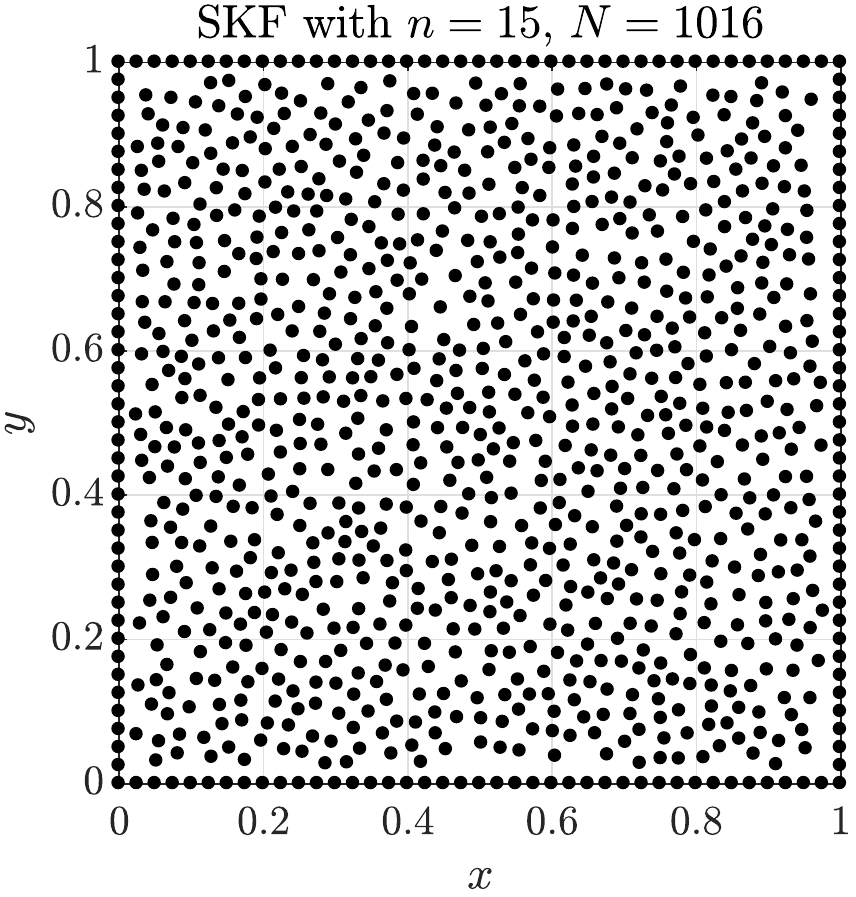}
  \includegraphics[width=0.24\linewidth]{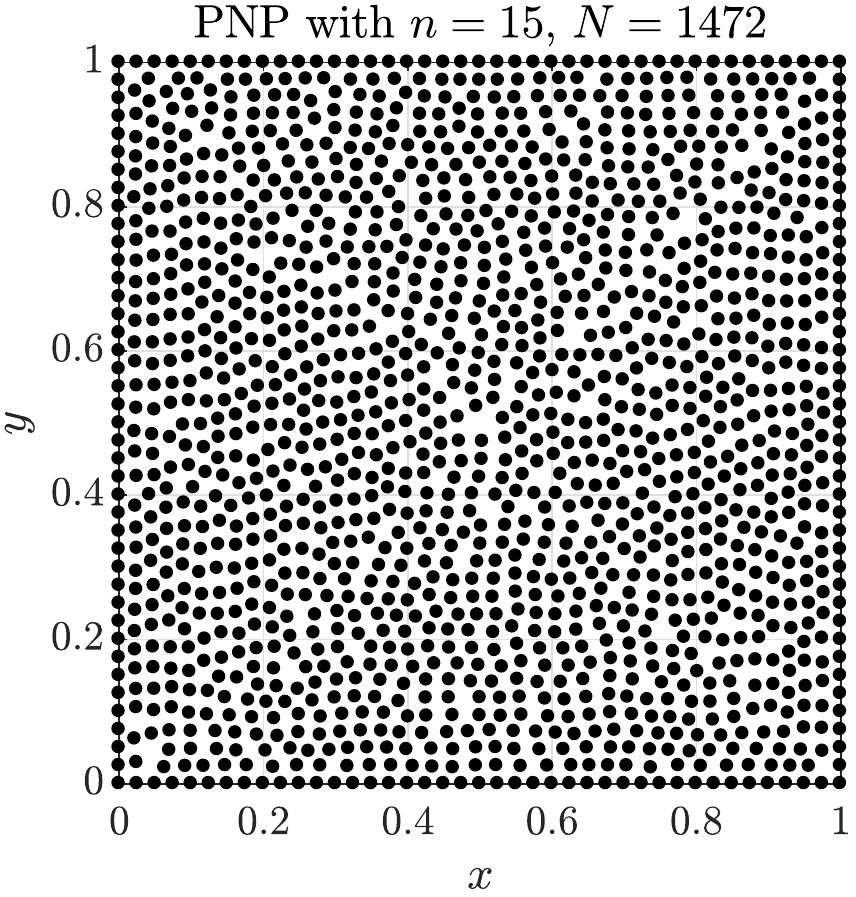}
  \includegraphics[width=0.24\linewidth]{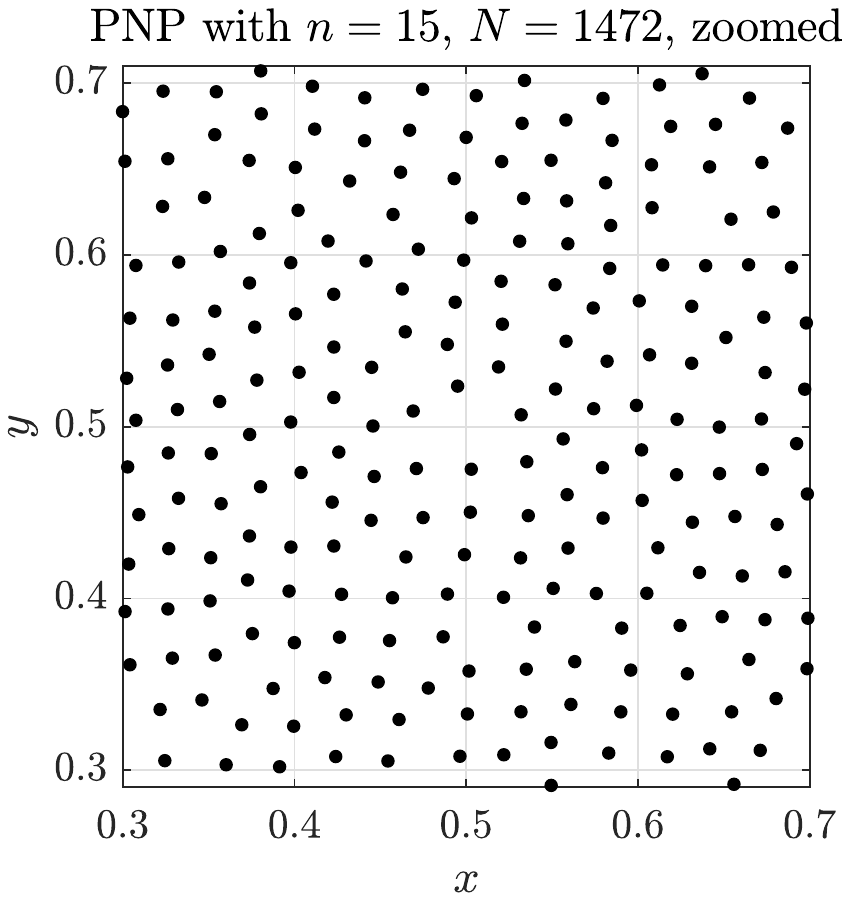}
  \caption{Node distributions on the unit square $[0, 1]^2$ with $h=0.025$
    generated with different algorithms. Rightmost figure shows the enlarged PNP
    distribution in the center where the advancing fronts meet.}
  \label{fig:sample-square}
\end{figure}

We begin our analyses by comparing the three algorithms on the unit square $[0, 1] \times [0, 1]$.
Node distributions were generated using constant density $h = 0.025$ and the expected number of
nodes is $N(h) = 1600$. Node distribution for all three algorithms are shown in~\cref{fig:sample-square}.
Parameters $n$ for various algorithms were chosen as recommended in their respective papers ($n=5$ for FF
and $n=15$ for SKF), with $n=15$ also being used for the algorithm presented in this paper.

SKF algorithm generated substantially less nodes than the other two. It also has
significant gaps between the boundary and internal nodes as well as visually more
irregular distributions. FF algorithm generates a smooth distribution without any
significant defects in the interior.
PNP algorithm exhibits gaps on diagonals, where advancing fronts from the sides have merged,
but behaves better near the boundaries. The part of the distribution where the advancing fronts
meet is shown in rightmost panel in Figure~\ref{fig:sample-square} to give a better perspective on
the size of the gaps.

In terms of the number of nodes, FF gives the best
result, since it produced only $45$ nodes less than expected, followed by PNP that produces
$128$ less nodes. The worst performance is demonstrated by SKF with deficiency of $573$ nodes.

\begin{figure}[h]
  \centering
  \includegraphics[width=0.32\linewidth]{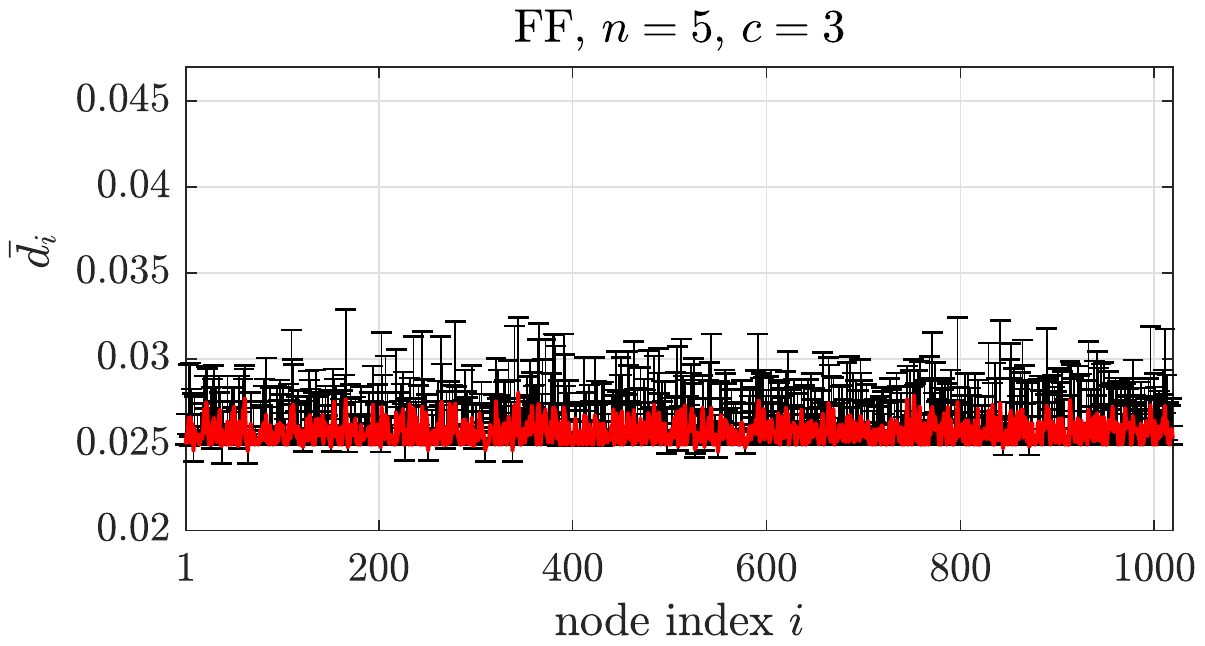}
  \includegraphics[width=0.32\linewidth]{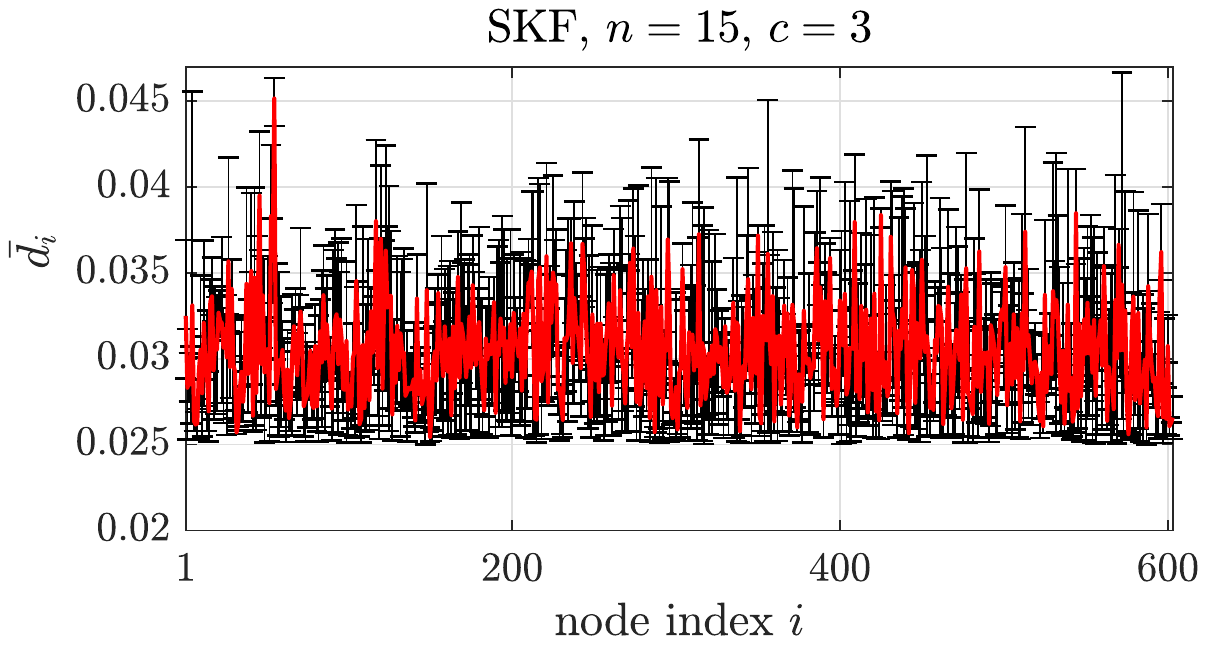}
  \includegraphics[width=0.32\linewidth]{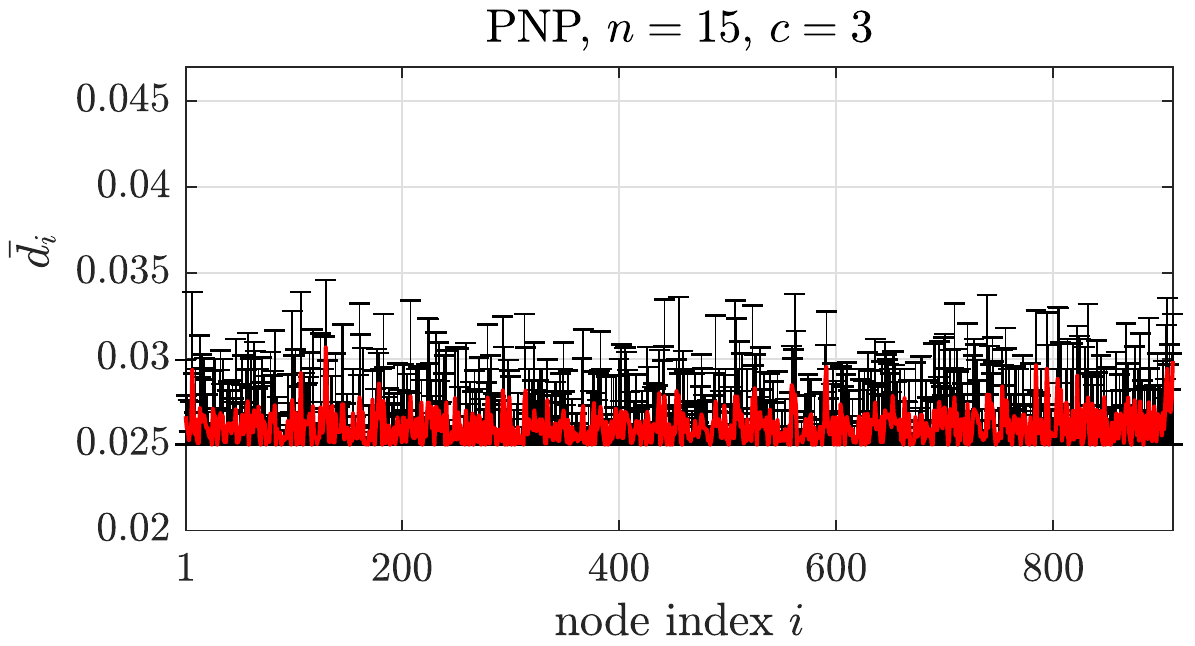}
  \caption{Average distances to the nearest neighbors for internal nodes.
    Error bars show minimal and maximal distances to three nearest neighbors.}
  \label{fig:support-dist}
\end{figure}

To analyze local regularity, distances to nearest neighbors are observed in the interior of the
domain. For each node $p_i$ at least $2h$ away from the boundary,
its $c$ closest neighbors (excluding $i$ itself) are found and denoted by $p_{i,j}$
for $j = 1, \ldots, c$ with distances to these neighbors computed as
$d_{i,j} = \|p_i - p_{i,j}\|$.
\Cref{fig:support-dist} shows average distances of each node to its three
closest neighbors, i.e.\ the plot of $\bar{d}_i = \frac13 \sum_{j=1}^3 d_{i,j}$
for each considered node $p_i$. Along with the average distance, the interval
$[d_i^{\text{min}}, d_i^{\text{max}}]$ is shown, where
\[
d_i^{\text{min}} = \min_{j=1,2,3} d_{i,j}, \qquad
d_i^{\text{max}} = \max_{j=1,2,3} d_{i,j}.
\]

FF and PNP algorithms show similar behavior, with average distance being close to $h$
with little variability between distances to closest neighbors.
FF algorithm has a few nodes a bit closer than $h$ together, but keeps the
internodal distance closer to $h$ and with less spread than PNP.
SKF algorithm performs worse with most of its distances to $c$ nearest neighbors
being on average closer to $0.03$ and with a significantly larger spread.
The numerical results representing these quantities are shown in~\cref{tab:supp-dist}.
The first two columns of the table demonstrate that the prescribed nodal spacing $h$ is much
better obeyed in PNP and FF algorithms, and the last column shows that
the average spread of the internodal distances in SKF algorithm is more than two
times greater than in FF and PNP.

Besides distances to the nearest neighbors, we can also take a look at the empty space
between the generated nodes. This can be done by computing the Voronoi diagram vertices
$v_j$ that lie inside the domain and observing the diameters $s_j$
of the largest circles centered at $v_j$ not containing any nodes. Formally, $s_j$ are given as
\begin{equation}
  s_j = 2 \min_i \|v_j - p_i\|.
\end{equation}
Note that the largest value of $s_j$ is the diameter of the largest empty circle. The basic
statistics os $s_j$ for the three considered algorithms are presented
in~\cref{tab:supp-dist}.

\begin{table}[h]
  \renewcommand{\arraystretch}{1.2}
  \centering
  \caption{Numerical quantities related to internodal distance and hole regularity.}
  \label{tab:supp-dist}
  \begin{tabular}{c|c|c|c|c|c|c}
    alg. & $\operatorname{mean}\bar{d}_i$ & $\operatorname{std}\bar{d}_i$ &
    $\operatorname{mean}(d_i^{\text{max}} - d_i^{\text{min}})$ & $\min s_j$ &
    $\operatorname{mean}s_j$ & $\max s_j$ \\ \hline
    FF  & 0.02575 & 0.00065 & 0.00208 & 0.028071 & 0.03438 & 0.04352 \\
    SKF & 0.03042 & 0.00275 & 0.02894 & 0.029737 & 0.04470 & 0.07008 \\
    PNP & 0.02604 & 0.00086 & 0.00276 & 0.028949 & 0.03568 & 0.05164 \\
  \end{tabular}
\end{table}

\begin{figure}[h]
  \centering
  \includegraphics[width=0.26\linewidth]{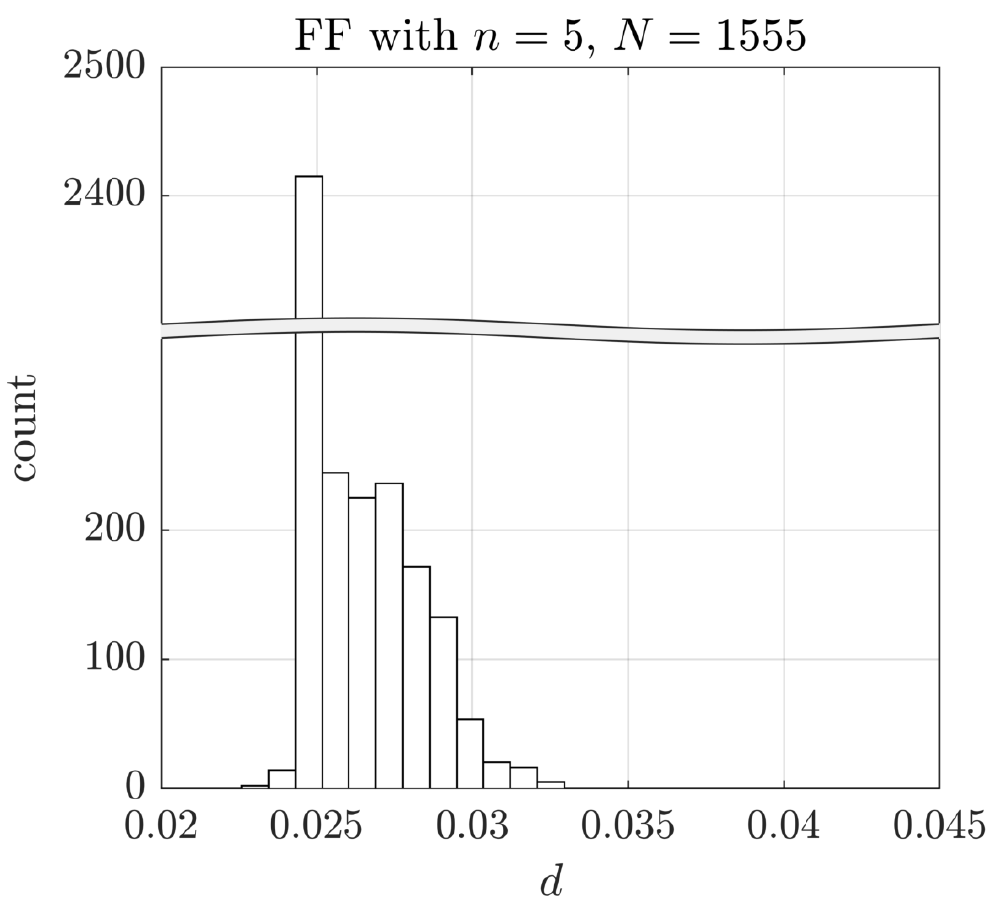}
  \includegraphics[width=0.26\linewidth]{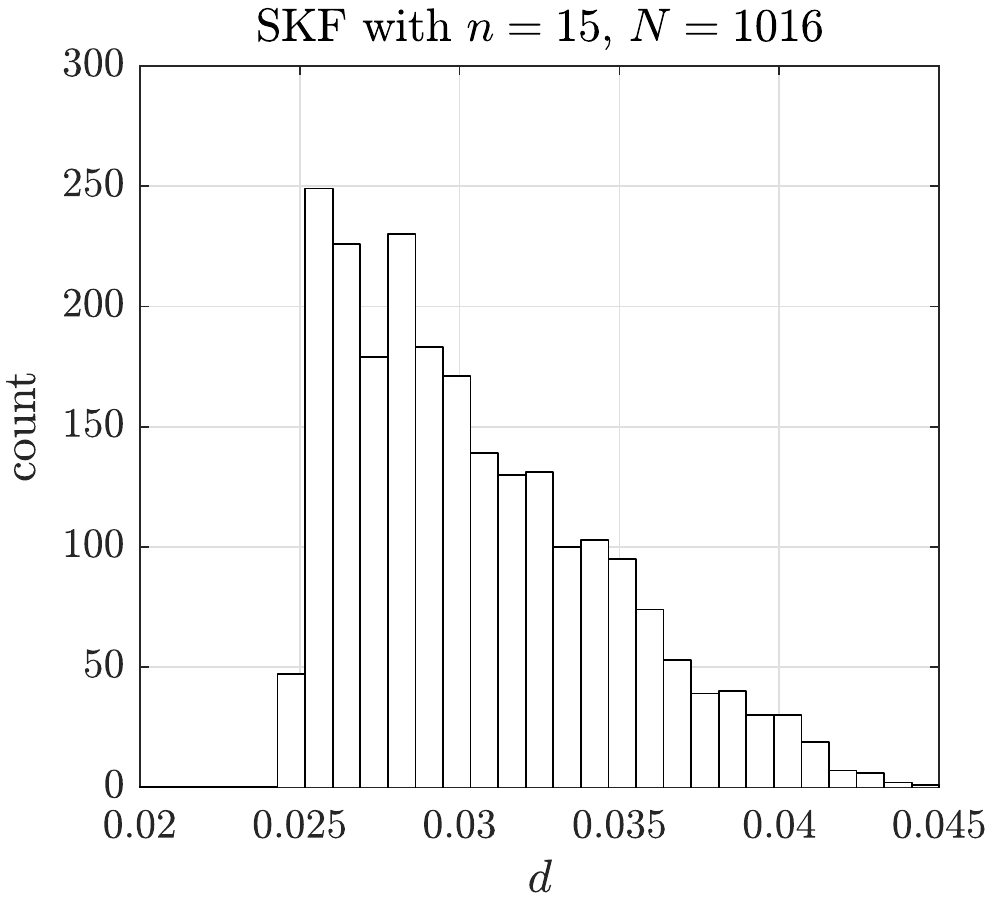}
  \includegraphics[width=0.26\linewidth]{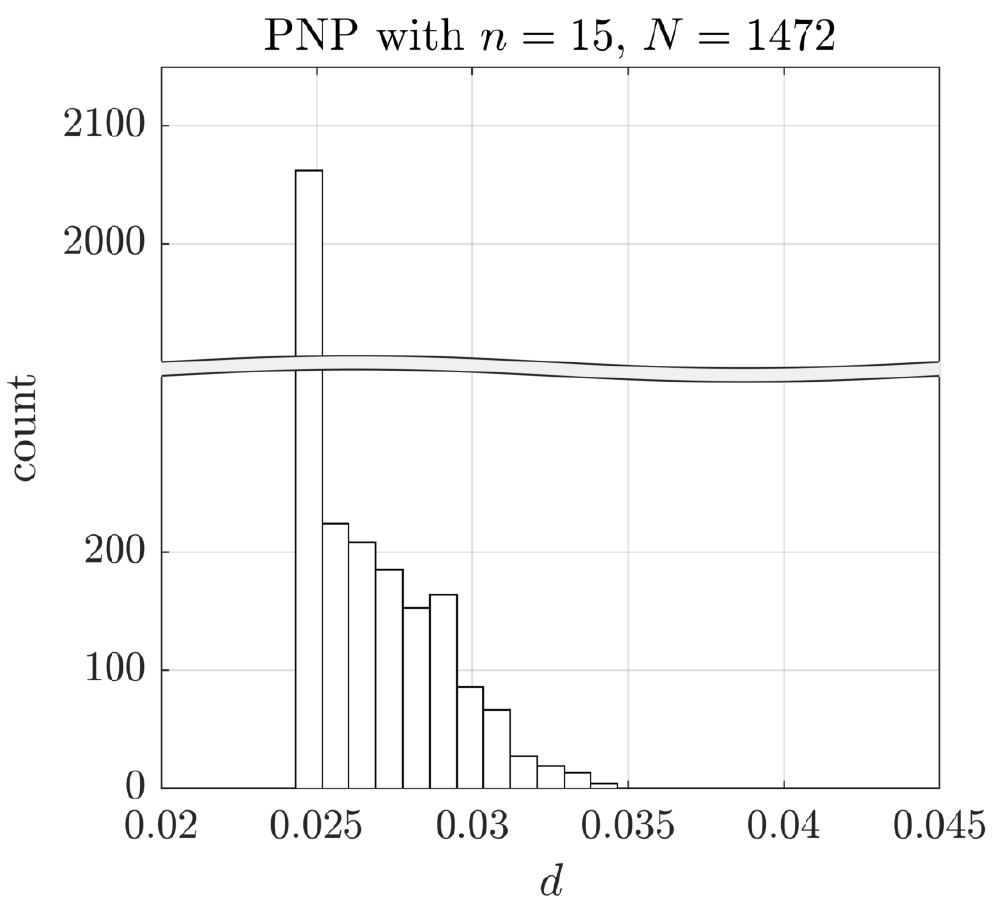}
  \caption{Histogram of distances to three nearest neighbors for node distributions
    on unit square $[0, 1]^2$ with $h=0.025$.}
  \label{fig:sample-square-hists}
\end{figure}

Additional insight is offered with histograms of distances to three nearest neighbors
(\cref{fig:sample-square-hists}). As expected, the largest count is in the bin around $h$. PNP and
SKF algorithms have no distances in bins below $h$, however the FF algorithm does put a small
number of nodes at a distance less than $h$ (see~\cref{sec:min-sp-req}). The irregularities
visible in the SKF algorithm distribution in~\cref{fig:sample-square} are reflected in the
histogram. The histogram has a much heavier tail than PNP and FF histograms, with far less nodes
exactly at distance $h$. PNP and FF histograms show more tightly packed distributions with slimmer
tails, however  the tail of PNP histogram is a bit longer and more spread out.

\begin{figure}[h!]
  \centering
  \includegraphics[width=0.25\linewidth]{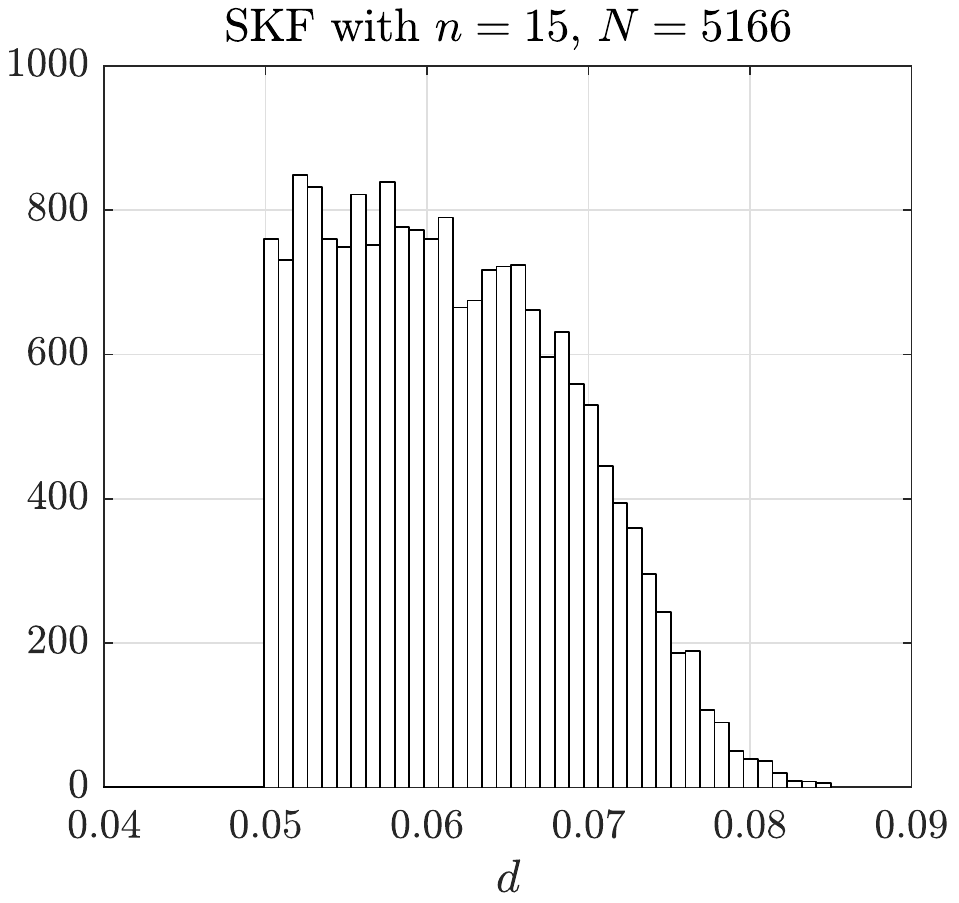}
  \includegraphics[width=0.25\linewidth]{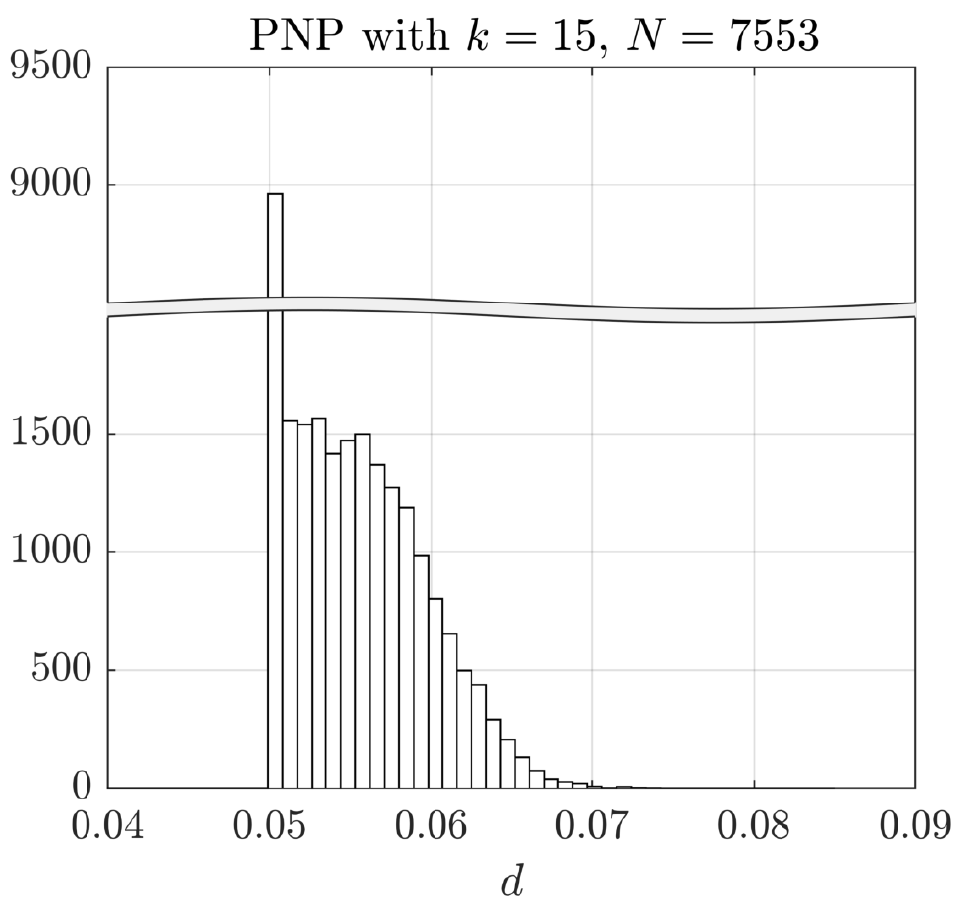}
  \caption{Histogram of distances to six nearest neighbors for internal nodes
  of distributions on unit cube $[0, 1]^3$ with $h=0.05$.}
  \label{fig:sample-square-3d-hists}
\end{figure}

Next, the PNP and SKF algorithms are compared in three dimensions.
The unit cube $[0, 1] \times [0, 1] \times [0, 1]$ is filled
with a constant density $h = 0.05$, starting from the boundary in the PNP case.
The expected number of nodes is $N(h) = 8000$. Histograms of distances to the closest $c=6$
nodes for internal nodes are shown in~\cref{fig:sample-square-3d-hists} for PNP and SKF
algorithms.

The histograms behave similarly to their 2-D counterparts. SKF algorithm again generated
significantly less nodes than the PNP algorithm. PNP has a large number of neighbors
at distance $h$ and a lighter tail, while the distances in SKF case are more spread out.

Further visual confirmation of regularity for variable density cases is
demonstrated in~\cref{sec:variable} (see~\cref{fig:variable},
\cref{fig:contact} and \cref{fig:sample-density}),
and more importantly, by the solutions of PDEs on generated node
sets~\cite{Fornberg_Flyer_2015,shankar2018robust,slak2018refined},
thus confirming sufficient local regularity.
Additionally, \cref{sec:num} considers sample solutions to PDE examples and
discusses accuracy, eigenvalue stability and convergence properties. Our experiments have
shown that SKF distributions cause stability problems when using small
stencils, such as e.g.\ closest 7 nodes. The likely cause of this instability is higher node
irregularity in SKF node distributions. PNP and FF distributions had no problems with
small stencils.

\subsection{Minimal spacing requirements}
\label{sec:min-sp-req}
Point~\ref{itm:min-dx} discusses minimal spacing guarantees. Provable minimal spacing guarantees
are very desirable, since nodes that are positioned too closely can effect the stability
of strong form methods. FF algorithm does not strictly respect the spacing $h$.
When running the algorithm with $h=0.005$ on a unit square $[0, 1]^2$,
some pairs of points in the domain interior were closer than $0.95 h$.
Although the violations do not appear to be significant and do not affect the
quality in practice, no bound of form $\|p_j-p_i\| \geq \alpha h$, for $\alpha > 0$ and
$i \neq j$ is known. 

%
%

SKF algorithm enforces the spacing between nodes to be greater than or equal to $h$ in the interior
and on the boundary, both times leveraging specialized spatial search structures.
The algorithm thus has the usual minimal spacing guarantee for constant nodal spacing:
\begin{equation}
\|p-q\| \geq h
\end{equation}
for $p\neq q$.

Similar argument can be made for PNP algorithm: each new candidate is checked using a $k$-d
tree against all previous ones, proving the minimal spacing guarantee for constant $h$.
For variable $h$, the above argument yields the bound
\begin{equation}
\|p_i - p_j\| \geq \min_{p \in \Omega} h(p)
\end{equation}
for $i \neq j$. This bound is dependent on a global property of $h$ and can be very coarse.
More precise, local bounds when considering spatially variable distributions are defined by
Mitchell et al.~\cite{mitchell2012variable}. If an ordered list of points, numbered 1 to $N$, is
considered, then the minimal spacing guarantee, called the \emph{empty disk property}, is
satisfied if
\begin{equation}
\|p_i - p_j\| \geq f(p_i, p_j),
\label{eq:empty-disk}
\end{equation}
for $1 \leq i < j \leq N$, where $f$ is a function evaluated at previously accepted
node $p_i$ and new candidate $p_j$. Four basic variations were proposed,
based on which point's spacing is taken into account when positioning new candidates:
\begin{itemize}
  \item \emph{Prior-disks:} $f(p_i, p_j) = h(p_i)$,
  \item \emph{Current-disks:} $f(p_i, p_j) = h(p_j)$,
  \item \emph{Bigger-disks:} $f(p_i, p_j) = \max\{h(p_i), h(p_j)\}$,
  \item \emph{Smaller-disks:} $f(p_i, p_j) = \min\{h(p_i), h(p_j)\}$.
\end{itemize}

The PNP procedure satisfies neither of this variations. The following proposition
establishes a version of the empty disk property~\eqref{eq:empty-disk} of PNP.
\begin{preposition}
  Let the points $p_i$, $i = 1, \ldots, N\!$, be a list of nodes generated by~\cref{alg:pds},
  where first $N_b$ nodes were given as initial nodes.
  The minimal spacing inequality
  \begin{equation}
  \|p_k - p_j\| \geq h(p_{\beta(j)})
  \end{equation}
  holds for all $N_b \leq k < j < N$. The function $\beta$ represents the predecessor function.
\end{preposition}
\begin{proof}
  \Cref{alg:pds} begins with $N_b$ initial nodes. Each candidate
  is generated from a unique existing node, thus giving rise to a predecessor-successor relation.
  Predecessor function $\beta\colon \{N_b+1, \ldots, N\} \to \{1, \ldots N\}$ for an accepted candidate
  $p_j$ that was generated from $p_i$ is defined as $\beta(j) = i$. Note that predecessors for
  the first $N_b$ initially given points are not defined.

  Consider an
  accepted candidate $p_j$, generated from a node $p_i$. The candidate was generated
  at a distance $h(p_i)$ from $p_i$, thus satisfying the equality
  \begin{equation}
  \|p_i - p_j\| = h(p_i) = h(p_{\beta(j)}).
  \end{equation}
  In particular, this means that~\cref{alg:pds} satisfies the \emph{prior-disks}
  property for predecessor-successor pairs.
  The distance $d$ to the nearest neighbor of $p_j$ among already accepted nodes is then found
  and if $d \geq h(p_i)$, the candidate is accepted. This means that the following inequality holds
  for all $k < j$:
  \begin{equation}
  \|p_k - p_j\| \geq d \geq h(p_i) = h(p_{\beta(j)}),
  \end{equation}
  establishing the desired property.
\end{proof}

\subsection{Spatial variability}
\label{sec:variable}
\begin{figure}[h]
  \centering
  \includegraphics[width=0.18\linewidth]{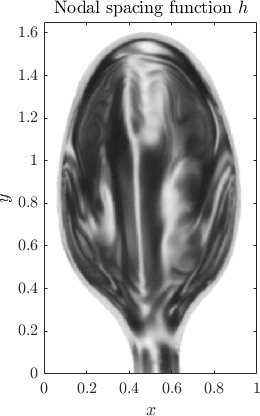}
  \includegraphics[width=0.18\linewidth]{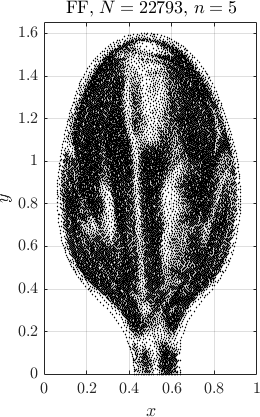}
  \includegraphics[width=0.18\linewidth]{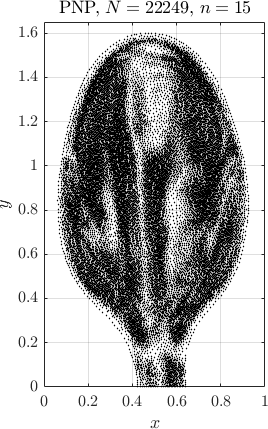} \\[1ex]
  \includegraphics[width=0.18\linewidth]{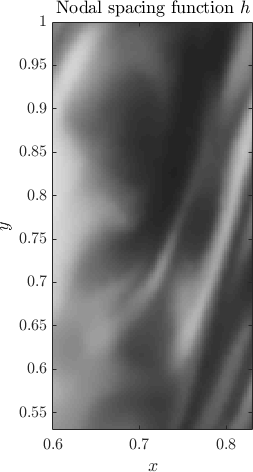}
  \includegraphics[width=0.18\linewidth]{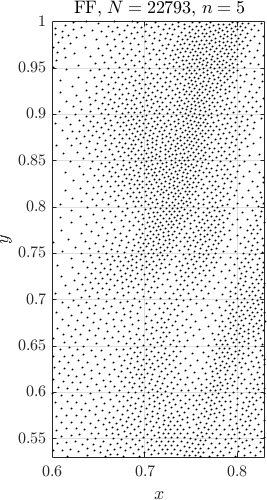}
  \includegraphics[width=0.18\linewidth]{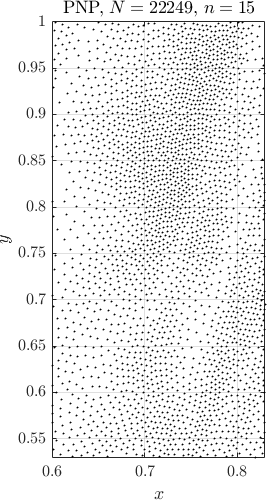}
  \caption{Illustration of variable density node sampling, with the nodal spacing
     function $h$ obtained from the image on the left using~\eqref{eq:img-to-h}.
     Enlarged variants are present to better asses the node quality.}
  \label{fig:sample-density}
\end{figure}

An important feature of FF and PNP algorithms is the ability to generate node sets with variable
nodal spacing on irregular domains.  SKF algorithm does not support variable nodal spacing
and is excluded from this analysis.
As an example, the image shown in top left corner of~\cref{fig:sample-density} is chosen
as a source for the nodal spacing function $h$.
The image is a modified version of an image showing stress distribution in a plastic spoon under
a photoelasticity experiment~\cite{spoon}. It features an irregular domain and rapidly varying dark
and light regions, which presents a more challenging case usually found in PDE discretizations.
The conversion from gray levels to the nodal spacing function is the same as used by Fornberg and
Flyer~\cite{fornberg2015fast}. Normalization factor $h_0 = 1.5$ was used to adjust the number of nodes
for maximal visibility.
The nodal spacing function $h$ is thus constructed from the image as
\begin{equation}
  \label{eq:img-to-h}
h(x, y) = h_0\, s\left(\frac{I_{\lfloor wx\rfloor, \lfloor wy \rfloor}}{255}\right), \; s(g) =
0.002+0.006\,g+0.012\, g^8,
\end{equation}
where $I_{ij}$ represents the grey level, ranging from 0 to 255, of the pixel in
the $i$-th row and the $j$-th column of the image and $w$ is the width of the image.
The node distributions obtained by filling the spoon shape with
aforementioned density using PNP and FF algorithms are shown in the first row
of~\cref{fig:sample-density}. The bottom row shows an enlarged portion of the
image and the corresponding distributions, so that individual nodes are visible
for easier visual assessment.

Both generated node sets conform to the supplied nodal spacing function.
The total number of nodes is similar in both cases, with PNP having
fewer nodes than FF. Enlarged portions show that PNP and FF distributions
are locally regular, visually similar and respect the variable nodal spacing
function $h$.

\begin{figure}[h!]
  \centering
  \includegraphics[width=0.58\linewidth]{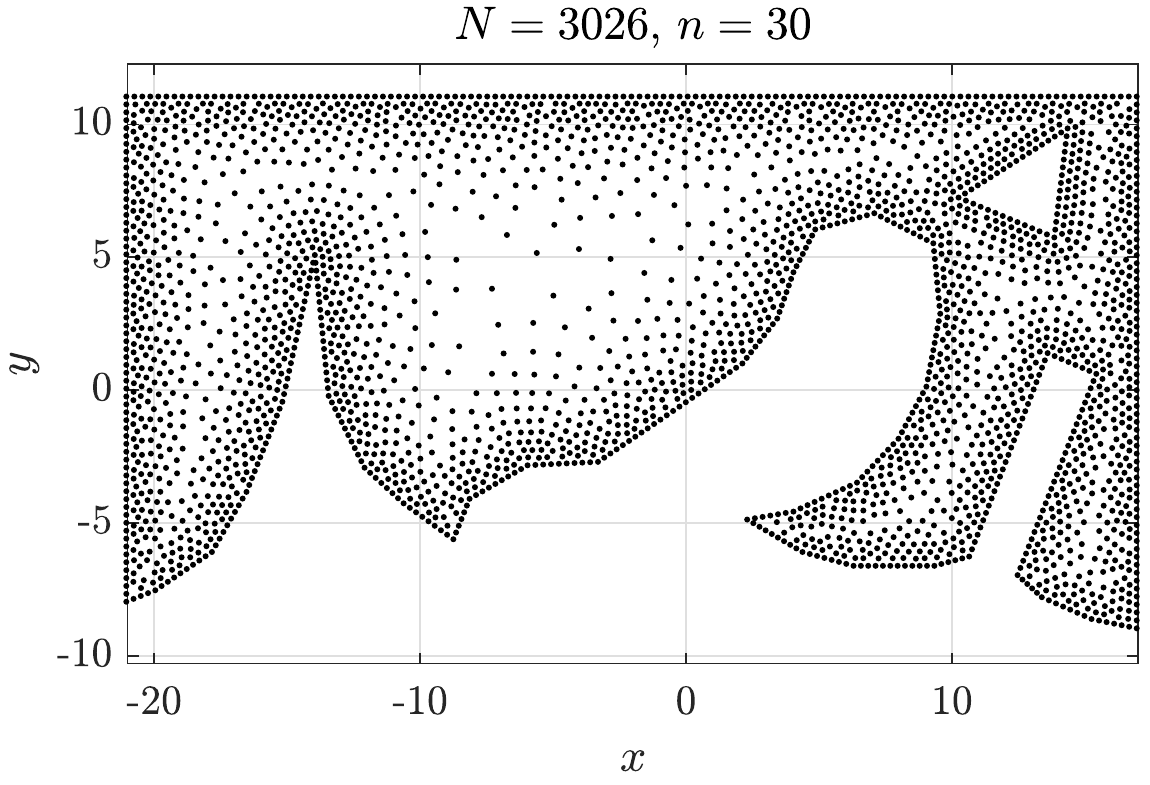}
  \includegraphics[width=0.4\linewidth]{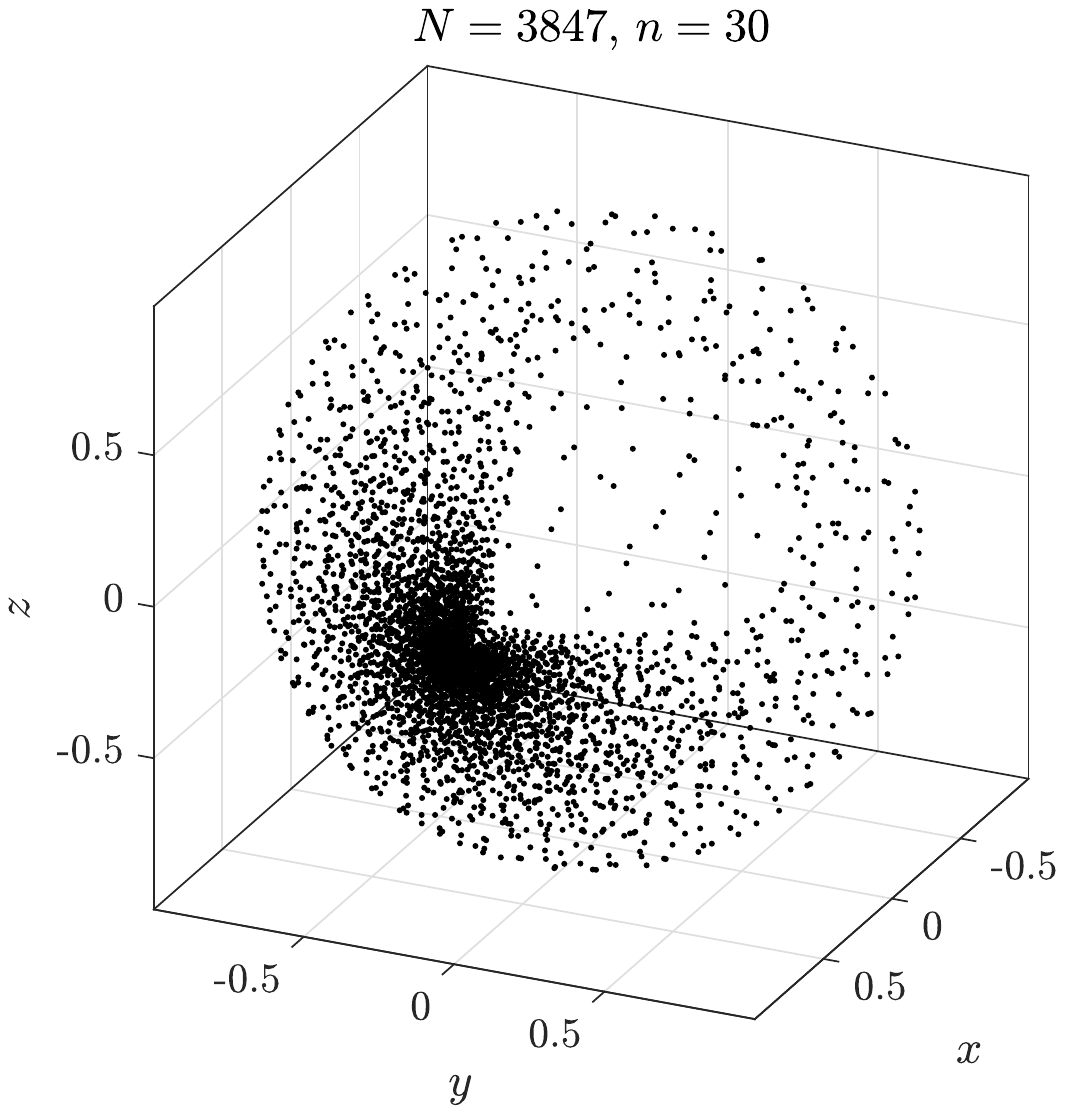}
  \caption{Example of generated variable density distributions for non-convex
  domain with non-trivial boundaries. }
  \label{fig:variable}
\end{figure}

\begin{figure}[h!]
  \centering
  \includegraphics[width=0.4\linewidth]{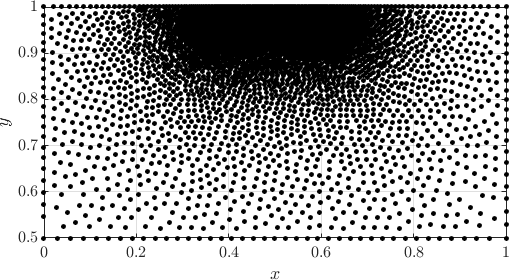}
  \includegraphics[width=0.4\linewidth]{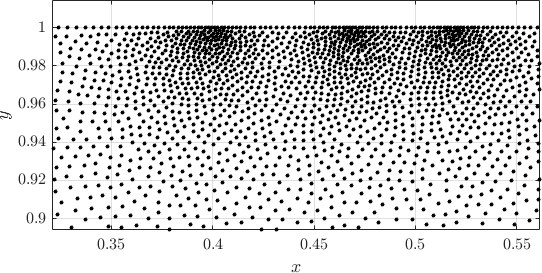}
  \includegraphics[width=0.4\linewidth]{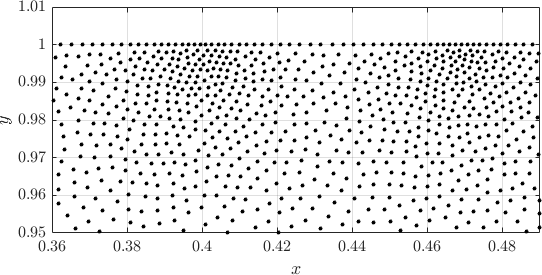}
  \includegraphics[width=0.4\linewidth]{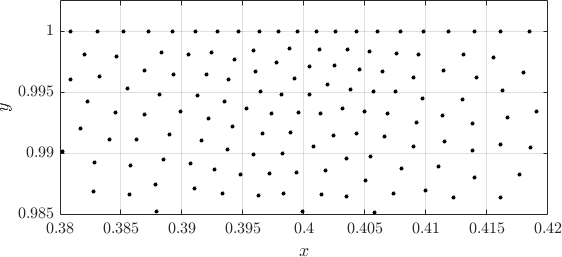}
  \caption{Discretisation of a contact region and successive enlargements.}
  \label{fig:contact}
\end{figure}

Further examples of 2D and 3D spatially variable distributions are shown
in~\cref{fig:variable}.
The 2D domain is a non-convex polygon with a hole and the 3D domain is a spherical shell
with one of the octants cut out.
\Cref{fig:contact} displays successive enlargement of a
nodal distribution used to solve a contact problem~\cite{slak2018adaptive}, which
illustrates the graded nature of the refinement and its local regularity in the most
zoomed panel.

\subsection{Computational efficiency and scalability}
\label{sec:compare-comp-eff}
Point~\ref{itm:time} concerns computational efficiency in two aspects: theoretical
time complexity and execution time.

The time complexity of the FF algorithm is proven in~\cref{sec:ff-time}
and given by~\eqref{eq:ff-time}
\begin{equation}
T_{\text{FF}} = O\left(n\left(\frac{|\bb{\Omega}|}{|\Omega|}N\right)^{1.5}\right)
\end{equation}
for constant spacing $h$ and is similar for variable spacing.
There are no immediate benefits if $h$ is assumed to be constant.
The SKF algorithm benefits from the assumption of constant $h$ and has time complexity
given by~\eqref{eq:skf-time} in~\cref{sec:skf-time}:
\begin{equation}
T_{\text{SKF}} = O\left(\frac{\text{obb}(\Omega)|}{|\Omega|} nN \right).
\end{equation}

PNP algorithm has time complexity
\begin{equation}
T_{\text{PNP}} = O(nN\log N),
\end{equation}
as analyzed in~\cref{sec:pds-time-complexity}. If $h$ is assumed constant,
the time complexity is further reduced to $O(\frac{|\text{obb}(\Omega)|}{|\Omega|}N + nN)$
using grid spatial search structure and even to $O(Nn)$ using hashing for irregular domains.

PNP algorithm is better for a domain irregularity factor compared to
both SKF and FF algorithms.
In case of constant $h$ it shares the same remaining factor $Nn$ with SKF
and for variable densities it is strictly better than FF.

Next, we compare the running time and scalability of proposed algorithms.
All time measurements were done on a laptop computer with an
\texttt{Intel(R) Core(TM) i7-7700HQ CPU @ 2.80GHz}
processor and 16 GB DDR4 RAM. Code was compiled using \texttt{g++ (GCC) 8.1.1}
for Linux with \texttt{-std=c++11 -O3 -DNDEBUG} flags.

Note that we implemented all three algorithms in the same manner with great emphasis
on optimization of the code in order to provide a fair comparison.

\begin{figure}[h]
  \centering
  \includegraphics[width=0.49\linewidth]{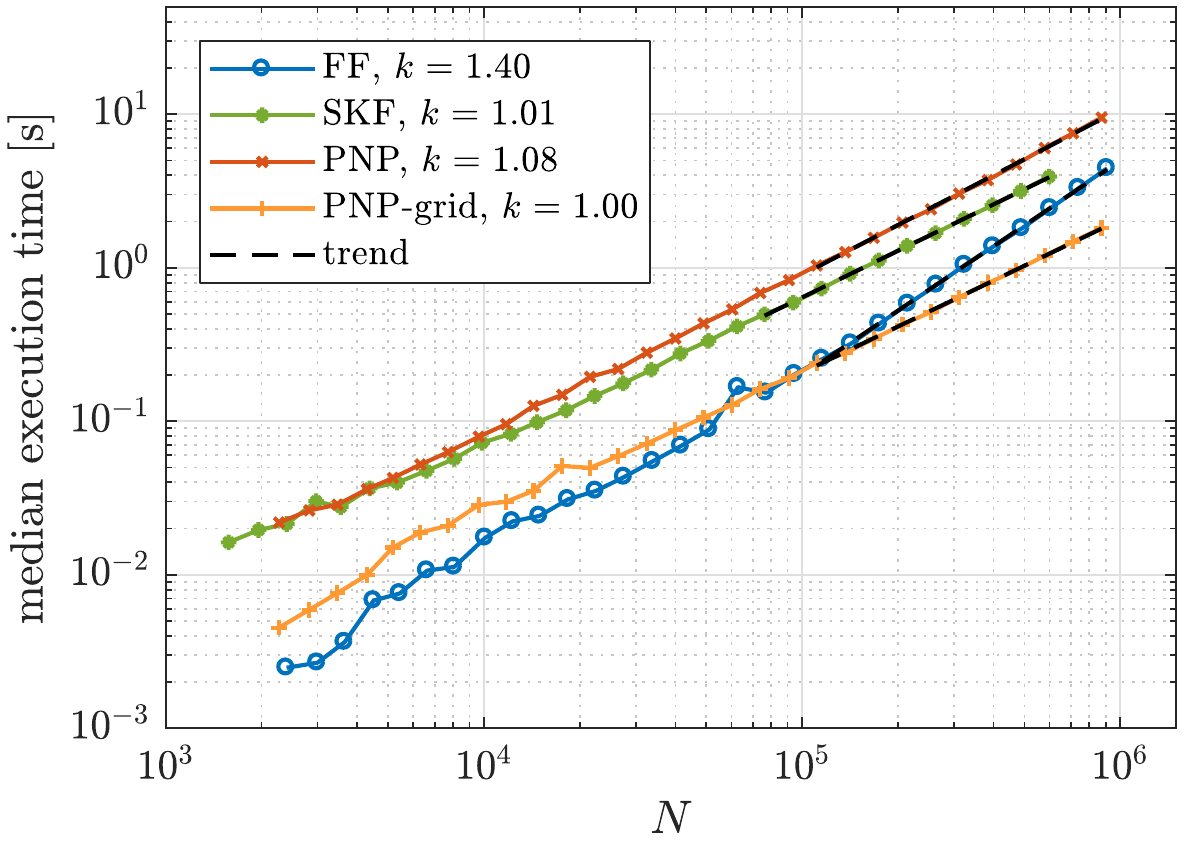}
  \caption{Execution times for the considered algorithms when filling $[0, 1]^2$ with
    successively smaller densities. Each data point represents a median of 10 runs. Standard deviation of  run times from the median was below 3\% in all cases. Value $k$ in the legend indicates slope of the line.}
  \label{fig:time}
\end{figure}

All algorithms were run on a unit square $[0, 1]^2$ with the
same parameters as in~\cref{sec:node-quality}. The nodal spacing function $h$
was varied as $h = \frac{1}{n}$, for such $n$ that the total number of nodes $N$ reached
approximately $N = 10^6$. Each run was executed 10 times and the median time was taken.
The results are shown in~\cref{fig:time}.

In 2-D the FF algorithm performs better than the others for small $N$. This is also expected,
as the algorithm generates nodes in a much simpler (and deterministic) way than the other two approaches.
SKF algorithm is next in terms of performance, with its grid-based search structure.
PNP algorithm is the slowest, due to the $k$-d tree search structure.
Nonetheless, $10^6$ nodes are generated in
5 to 10 seconds, which is significantly less than the time that would be spent on solving
the PDE on these nodes.

The trends for large $N$ coincide with the theoretical time complexities with SKF being an
$O(N)$ algorithm, PNP being an $O(N \log N)$, and FF being $O(N\sqrt{N})$.

PNP was also run using the same grid search structure as SKF, denoted in~\cref{fig:time} by
``PNP-grid''. It shows a significant improvement over the use of $k$-d tree spatial search structure
and agrees with the predicted linear time complexity. This also shows that
PNP algorithm itself is about three times faster than SKF, when compared using the same search
structure and the same number of candidates. PNP with a gird-based structure also comes close to
FF for smaller $N$ and constant $h$.
Execution time of PNP and SKF algorithms was tested also in 3-D and the results are equivalent.

Additionally, we analyze the execution time of the three algorithms when dealing with irregular domains.
Both FF and SKF algorithms do no have
time complexity proportionate to $|\Omega|$, but rather to the volume of its (oriented)
bounding box, which can be arbitrarily larger. In practice this means that PNP algorithm
inherently benefits in execution time by a factor of $\frac{|\operatorname{bb}(\Omega)|}{|\Omega|}$.

\begin{figure}[h]
  \centering
  \includegraphics[width=0.4\linewidth]{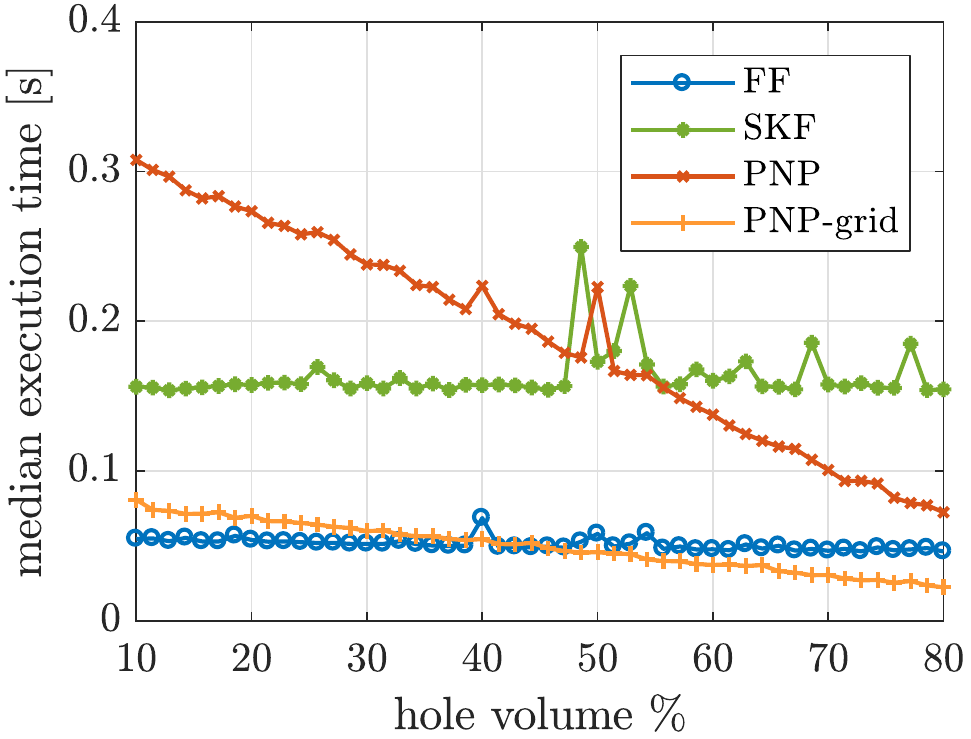}
  \caption{Execution time when filling domains $\Omega(\alpha)$ which have decreasing area.}
  \label{fig:shrinking-domain-results}
\end{figure}

This is illustrated in~\cref{fig:shrinking-domain-results}, which shows the execution time
of the considered algorithms when filling increasingly ``emptier'' domains
\begin{equation}
\Omega(\alpha) = [0, 1]^2 \; \setminus \; \Big(\frac{1}{2}-\alpha, \frac{1}{2}+\alpha\Big)^2.
\end{equation}

Domains $\Omega(\alpha)$ are chosen in such a way that the bounding box is equal to
$[0,1]^2$ for all $\alpha$ and that the limit of the ratio between the bounding box and domain
volume approaches zero as $\alpha$ approaches 1/2.

The difference in the behavior of execution time is substantial
and shows that both versions of PNP really scale with volume of $\Omega$,
while the execution time of FF and SKF remains almost constant, as predicted by time complexity
analysis. This means that around $30\,000$ nodes can be generated, for less than $8\,000$ to
remain in the final set.

\subsection{Compatibility with boundary discretizations}
The next point discusses the compatibility between interior and boundary discretizations.
All three algorithms treat boundary discretizations separately from discretizing the interior of $\Omega$.
Due to box-fill nature of FF, the generated discretization of the interior is chopped off at the boundary
of $\Omega$ when the boundary discretization is superimposed.
Nodes that are closer to the boundary nodes than a given threshold are discarded.
If the threshold is strictly $h$, gaps between the boundary and the interior discretization can occur.
The authors recommend setting the threshold to $h/2$ and preforming a few iterations of
a repel-type algorithm that is executed locally on the nodes near the boundary to smooth
the transition between both discretizations. SKF algorithm possesses a similar problem, but
deals with it differently. It generates internal nodes in a slightly reduced oriented bounding box,
which is computed from boundary nodes that were shifted by $h$ to the interior.
This prevents the generation of internal nodes too close to the boundary of the box;
however, the nodes still need to be tested for inclusion, which is done using
the shifted nodes and their normals. This causes gaps near the boundary, which can be
observed in the sample distribution in~\cref{fig:sample-square} (second panel).

PNP algorithm bypasses the aforementioned problems altogether, by offering the option to
use the boundary discretization as a starting point of the interior discretization and
thus allowing for a smooth transition near the boundary.
Similar irregularities to
those present near the boundary in SKF algorithm
are formed when the advancing fronts from the opposite sides meet,
but they appear in the interior of $\Omega$ (see~\cref{fig:sample-square}, rightmost panel),
where they have less impact on the stability of the solution. Consequently
no need to smooth the irregularities with expensive iterative repel techniques arose.

\subsection{Compatibility with irregular domains}
\label{sec:comp-irreg}

Another requirement deals with irregular domains. FF algorithm has a disadvantage
of being only able to fill axis-aligned boxes, which results in potentially a lot
of unnecessarily generated nodes. This approach is somewhat improved in the
SKF algorithm, where oriented bounding boxes are used, in general reducing the number of
generated nodes compared to FF.
The number of unnecessarily generated nodes could be reduced even further by decomposing an
unfavorably shaped domain into smaller domains, which can be better bounded by cuboids. The smaller
domains can then be filled separately and combined together, provided that the node generation
algorithm behaves well near boundaries. An appropriate domain decomposition would also enable
immediate parallel execution of the algorithm.

Of the three discussed algorithms only PNP never generates any unnecessary nodes in the
exterior of the given domain $\Omega$, never evaluates nodal spacing function
$h$ outside of $\Omega$ and has the property that the number total number of generated
nodes and the time complexity scale directly with $|\Omega|$.
The impact of unnecessary node generation outside $\Omega$ on the execution time
is illustrated in~\cref{sec:compare-comp-eff}; however, the slowdown introduced by bounding boxes
is in practice often acceptable.

The strength of the PNP algorithm which allows it to generate nodes only inside $\Omega$
can also become its disadvantage. If seed nodes are supplied only in one part of $\Omega$
and the domain has a bottleneck in the middle (such as an hourglass shape) of girth
approximately equal to nodal spacing $h$ in that area, the algorithm
might fail to advance through such bottleneck and would not generate any nodes in the other part.
The FF and SKF algorithms do not suffer from this problem, and it can also be circumvented
in PNP by supplying at least one seed node in each problematic part of the domain.

\subsection{Direction and dimension independence}
Points~\ref{itm:dim-ind} and~\ref{itm:dir-ind} deal with direction and dimension
independence. FF algorithm is only two-dimensional and directionally dependent,
because the advancing front progresses with respect to the increasing $y$ coordinate.
For inconveniently rotated or badly shaped domains, filling via increasing last coordinate
might perform badly. Choosing a filling direction is the first step of the algorithm,
and  it can have significant effect on the running time and the generated node distribution.
The algorithm is also not easily generalizable to higher dimensions,
as it is not immediately obvious how to extend the concept of the ``closest left'' point to
higher dimensional spaces.

The SKF algorithm is better in this aspect. Using PCA
it computes oriented bounded boxes, which provides independence from rotations.
The main parts of the SKF algorithm, i.e.\ PCA and
Poisson Disk Sampling, all work in arbitrary dimensions.
Similarly, all parts of the PNP algorithm are formulated for a general dimension $d$
and the formulation of the fill procedure is independent of the coordinate system.
The same is true for the implementation: there is a single implementation for
all values of $d$ and the space dimension can truly be a run-time parameter.
The coordinates of points are only accessed in the internals of the $k$-d tree
operations; all other expressions are coordinate-free.

\subsection{Free parameters}
Point~\ref{itm:free-param} states that the developed algorithm should aim to minimize the number of
free or tuning parameters. All three algorithms have a parameter influencing the number of
candidates, which represents a time-quality trade-off. Authors of FF set $n=5$ and
anything above has similar distributions with a higher execution time.
Authors of SKF analyze the effect of the number of candidates more precisely and
recommend $n = 15$ in 2-D and 3-D, with a higher number of candidates corresponding
to lower errors. For PNP algorithm we similarly recommend $n=15$ in 2-D, and $n=30$ in
3-D with increasing $n$ for higher dimensions. Anything above $n=30$ in 2-D
gives very similar results and is computationally wasteful.

\section{Solution of PDEs on generated nodes}
\label{sec:num}

\subsection{Poisson's equation}
The decisive factor of node distribution quality for strong form methods
is its ability to support construction of a good approximations of differential
operators.
A basic test of this ability is to solve the Poisson's equation on nodes
generated by all three algorithms and compare the accuracy of the solutions.

A $d$-dimensional boundary value problem
\begin{align}
\nabla^2 u &= f \quad \text{ in } \Omega = [0, 1]^d, \nonumber \\
u &= 0 \quad \text{ on } \partial \Omega \label{eq:lap}
\end{align}
with $u(x_1, \ldots, x_d) = \prod_{i=1}^d \sin(\pi x_i)$ and
$f(x_1, \ldots, x_d) = -d\pi^2\prod_{i=1}^d \sin(\pi x_i)$ is considered
in $d=2$ and $d=3$ dimensions.

The solution is obtained using the popular strong form RBF-FD
method~\cite{Fornberg_Flyer_2015, mavric, slak2018refined}.
Polyharmonic radial basis functions (PHS)
\begin{equation}
  \phi(r) =
  \begin{cases}
    r^k & k \text{ odd} \\
    r^k \log r & k \text{ even}
  \end{cases}
\end{equation}
with $k=3$ augmented with monomials up to order 2 are used to construct the
approximations
on a stencil of 15 closest nodes in 2-D
and 42 closest in 3-D. The final system is solved using BiCGSTAB iterative
algorithm
with tolerance $10^{-15}$ and 100 iterations with ILUT preconditioner with fill factor
20 and drop tolerance $10^{-5}$.

The $L^1$ error between the correct solution $u$ and obtained solution $u_h$
is evaluated on an independent uniform grid of points $G$,
three times denser than the densest
discretization used in solution of the problem, and computed as
\begin{equation}
  L^1 = \|u_h - u\|_1 \approx \frac{1}{|G|} \sum_{p \in G} |u(p) - u_h(p)|.
\end{equation}

\begin{figure}[h]
  \centering
  \includegraphics[width=0.35\linewidth]{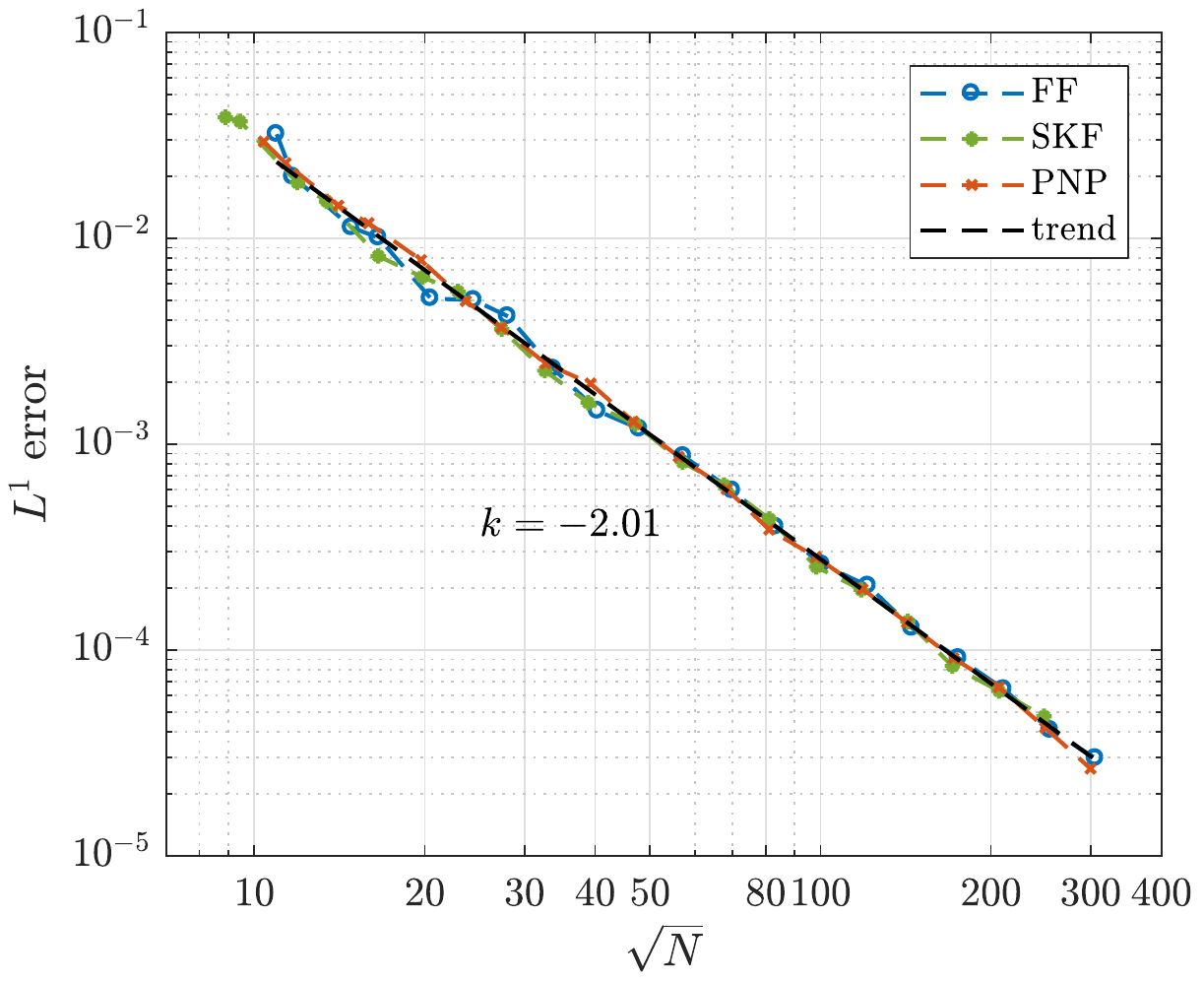}
  \includegraphics[width=0.35\linewidth]{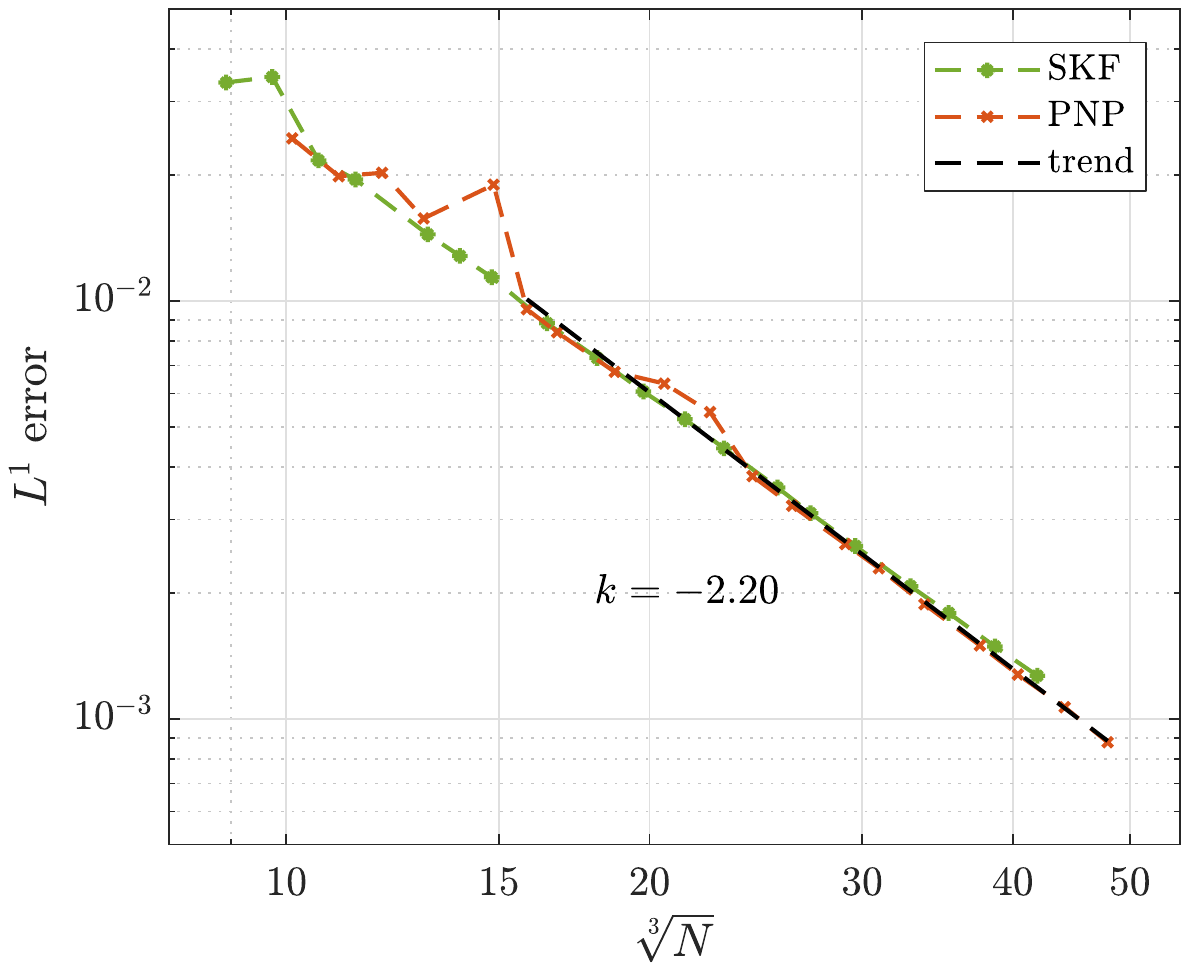}
  \caption{Accuracy of the numerical solution of~\eqref{eq:lap}
    for considered algorithms when filling $[0, 1]^d$ with successively smaller
    densities in 2D (left) and 3D (right).}
  \label{fig:sample-square-pde-acc}
\end{figure}

Node distributions generated by the three considered algorithms are tested
using the same
parameters as in~\cref{sec:node-quality} and~\cref{sec:compare-comp-eff}. The nodal spacing
function $h$ varies as $h = \frac{1}{n}$, for such $n$ that the total number of
nodes $N$
reached approximately $N = 10^5$. The results are shown in~\cref{fig:sample-square-pde-acc}.

In both 2D and 3D case we observe convergence with expected rates for large $N$.
All node sets give well-behaved solutions with very similar accuracy.
Similar results are obtained in 3D.
\subsection{Eigenvalue stability}
\begin{figure}[h]
  \centering
  \includegraphics[width=0.43\linewidth]{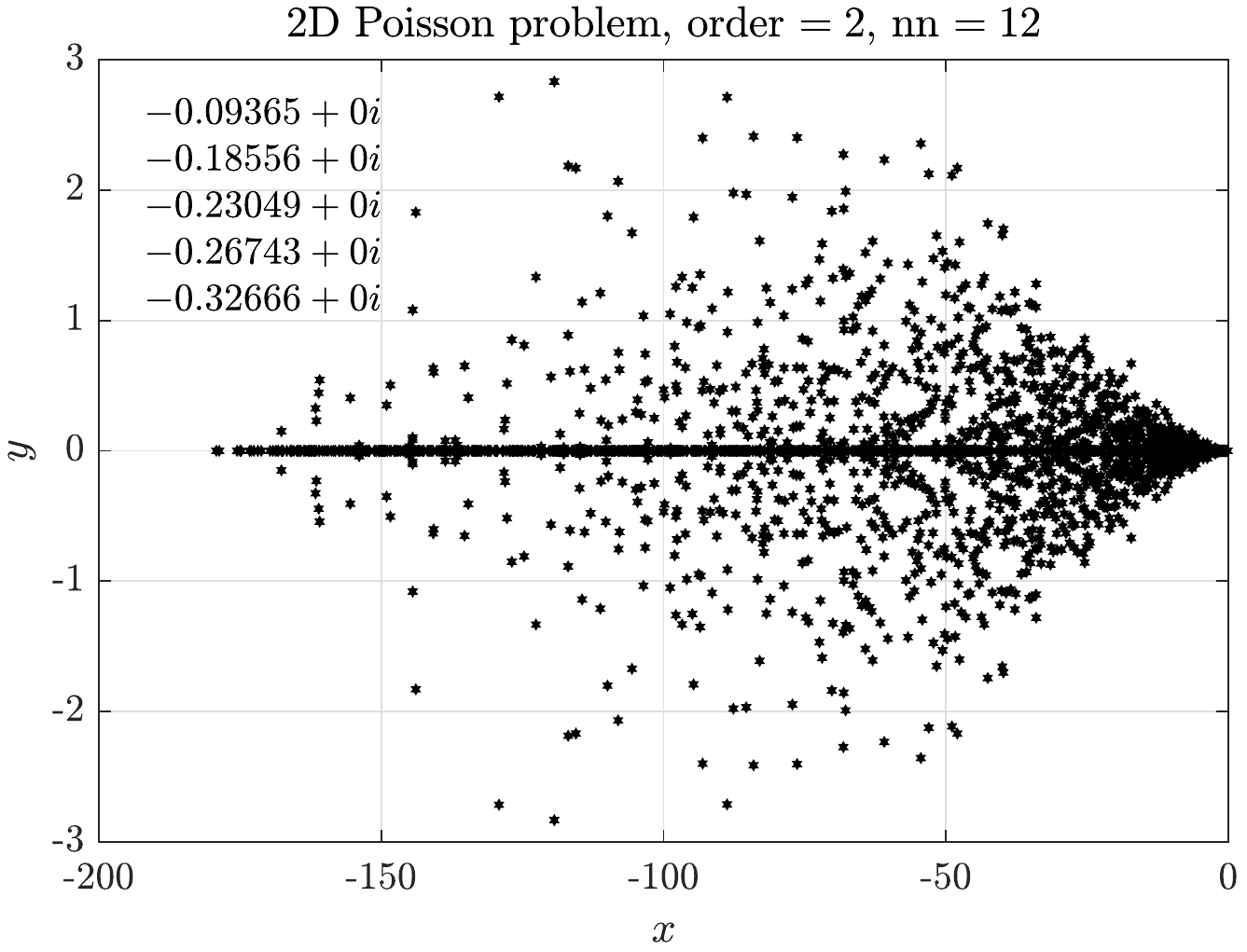}
  \includegraphics[width=0.45\linewidth]{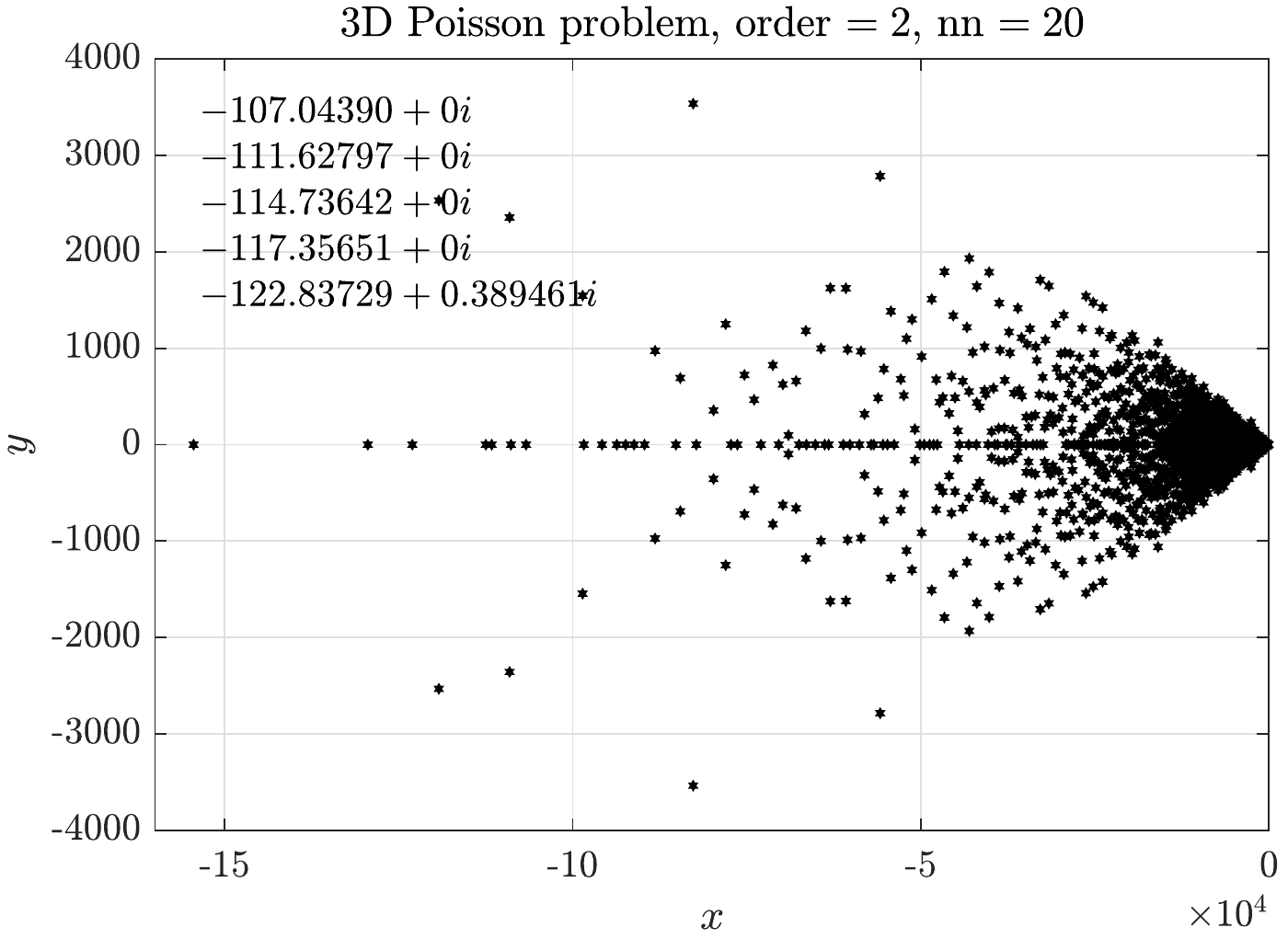}
  \caption{Spectra of the Laplacian operator discretized with RBF-FD PHS $r^3$,
  augmented with monomials or order 2 on PNP nodes. Note the different scales on the axes of both
  plots. Variable $\textup{nn}$ denotes the number of nearest neighbors used to construct the
  stencil. The 5 eigenvalues with the largest real parts are given in the top left corner of each
  plot.}
  \label{fig:eigen}
\end{figure}

An often observed property of numerical discretization methods is the
spectrum of discretized partial differential operator. For example, the spectrum of
discretized Laplace operator should have only eigenvalues with negative
real part, and a relatively small spread along the
imaginary axis~\cite{shankar2018robust}.
\Cref{fig:eigen} shows the spectrum of Laplace operator discretized with 2nd order
RBF-FD PHS on PNP nodes shown in~\cref{fig:variable}. There are no eigenvalues
with positive real part and also imaginary spread is relatively small, which
additionally confirms the stability of RBF-FD PHS differentiation on scattered nodes.

Additionally, we tested several different setups with different stencil sizes and approximation
orders on nodes distributed with all three positioning algorithms with minimal differences
observed in the spectrum.

\subsection{Thermo-fluid problem}

Finally, the PNP algorithm is tested on a more complex problem. The goal is to
demonstrate the capability of meshless solution procedure on PNP nodes solving
a transient non-linear convection dominated problem in 2D and 3D
irregular domain with mixed Dirichlet and Neumann boundary conditions. The
natural convection problem governed by coupled Navier-Stokes, mass continuity and heat
transfer equations is chosen for a test case. First, a well-known de Vahl Davis benchmark
test~\cite{De_Vahl_Davis_1983} is solved  to demonstrate correctness of the solution procedure,
both in 2D and 3D. Once we attain confidence in the solution procedure we extend the demonstration
to irregular domains.

The natural convection benchmark problem is governed by the following equations:
\begin{align}
\dpar{\b v}{t} + (\b v \cdot \nabla) \b v &= -\frac{1}{\rho} \nabla p + \frac{\mu}{\rho} \nabla^2
\b v + \frac{1}{\rho} \b b, \\
\nabla \cdot \b v &= 0, \\
\b b &= \rho (1-\beta(T - T_{\text{ref}})) \b g, \\
\dpar{T}{t} + \b v \cdot \nabla T  &= \frac{\lambda}{\rho c_p} \nabla^2 T,
\end{align}
where $\b v(u,v,w)$, $p$, $T$, $\mu$, $\lambda$, $c_p$, $\rho$, $\b g$,
$\beta$, $T_{\text{ref}}$ and
$\b b$ stand for velocity, pressure, temperature, viscosity, thermal conductivity,
specific heat, density, gravitational acceleration, coefficient of thermal expansion,
reference temperature for Boussinesq approximation, and body force,
respectively. The de Vahl Davis test is defined on a unit square domain $\Omega$,
where vertical walls are kept at constant temperatures
with $\Delta T$ difference between cold and hot side, while
horizontal walls are adiabatic. In generalization to 3D we assume also front
and back walls to be adiabatic~\cite{wang2017numerical}. No-slip velocity
boundary conditions are assumed on all walls. The problem is characterized
by Rayleigh (Ra) and Prandtl (Pr) numbers, defined as
\begin{equation}
  \textup{Pr} = \frac{\mu c_p}{\lambda},
  \textup{Ra} = \frac{g \beta \rho c_p \Delta T h^3}{\lambda \mu},
\end{equation}
with $h$ standing for characteristic length, in our case set to $1$. All cases
considered in this paper are computed at $\textup{Pr} = 0.71$.

The problem is solved with implicit time stepping, where each time step begins
with a computation of intermediate velocity ($\tilde{\b v}_2$)
\begin{equation}
  \tilde{\b v}_2 = \b v_1 + \Delta t \left[ -(\b v_1 \cdot \nabla) \tilde{\b v}_2 +
  \frac{\mu}{\rho} \nabla^2 \tilde{\b v}_2  + \frac{1}{\rho} \b b(T_1) \right].
\end{equation}
The computed velocity is coupled with mass continuity by an
iterative velocity-correction scheme, where it is assumed that the correction
depends only on the pressure term
\begin{equation}
\b v_2 = \tilde{\b v}_2 - \frac{\Delta t}{\rho} \nabla p.
\label{eq:PVI_1}
\end{equation}
Applying divergence on~\eqref{eq:PVI_1} yields a pressure Poisson equation
\begin{equation}
 \nabla^2 p = \frac{\rho}{\Delta t} \nabla \cdot \tilde{\b v}_2 \text{ in }
 \Omega,
 \quad \dpar{p}{n} = \frac{\rho}{\Delta t} \tilde{\b v}_2 \cdot \b n \text{ on
 } \partial \Omega,
 \quad \text{ subjected to } \int_{\Omega} p = 0,
 \label{eq:PVI_2}
\end{equation}
which is solved first to get the pressure field. With
computed pressure the
velocity is corrected following the~\eqref{eq:PVI_1}. Steps~\eqref{eq:PVI_2}
and~\eqref{eq:PVI_1} are repeated until the convergence
criterion is not met. Once the velocity is satisfactorily divergence free, the
temperature field, coupled with momentum equation through Boussinesq approximation,
is updated as
\begin{equation}
  T_2 = T_1 + \Delta t \left[ - \b v_2 \cdot \nabla T_2 + \frac{\lambda}{\rho c_p} \nabla^2 T_2
\right].
\end{equation}


All spatial operators are discretized using RBF-FD with $r^3$ PHS radial basis
functions, augmented with monomials up to order $2$, with the closest $25$
nodes used as a stencil. For the time discretization time step
$\Delta t=10^{-3}$ was used for all cases. Nodal distance $h=0.01$ is used for
simulations in 2D and $h=0.25$ for simulations in 3D. Boundaries with Neumann
boundary conditions are additionally treated with ghost
nodes~\cite{bayona2017augment}.

\begin{figure}[h]
  \centering
  \includegraphics[width=0.45\linewidth]{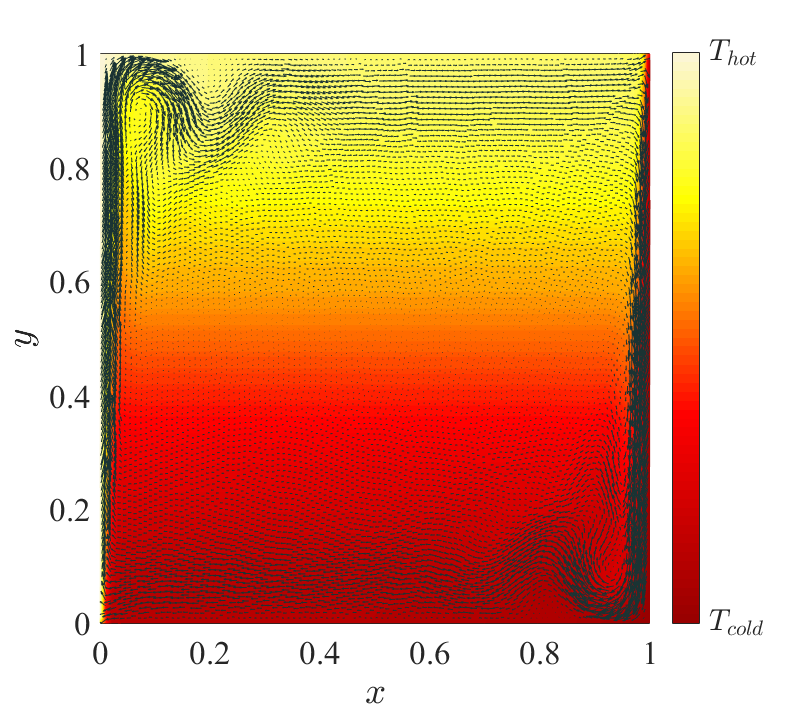}
  \raisebox{-12pt}{\includegraphics[width=0.48\linewidth]{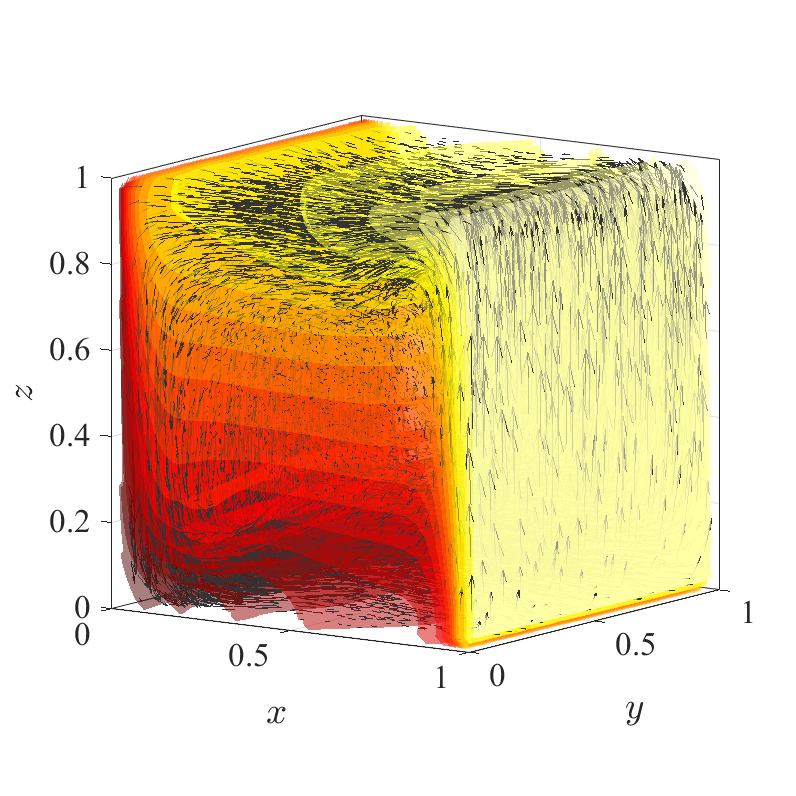}}
  \vspace{-12pt}
  \caption{Temperature contour and velocity quiver plots
    for $\textup{Ra} = 10^8$ case in 2D (left) and $\textup{Ra} =10^6$ case in 3D (right).}
  \label{fig:dvd-uni}
\end{figure}

In~\cref{fig:dvd-uni} steady state temperature contour and velocity quiver plots
for $\textup{Ra} = 10^8$ case in 2D and $\textup{Ra} = 10^6$ case in 3D are presented.
A more quantitative analysis is done by comparing characteristic values, i.e.\
peak positions and values of cross section velocities, with data available in
literature~\cite{couturier2000performance, kosec2008solution,
wang2017numerical, fusegi1991numerical}. We analyze six different cases,
namely $\textup{Ra} = 10^6,10^7,10^8$ in 2D, and $\textup{Ra} = 10^4,10^5,10^6$ in 3D.
The comparison in presented in~\cref{tab:ff-data}.

\begin{table}[h]
  \centering
  \caption{Comparison of results computed with RBF-FD on FF nodes and reference
  data. }
  \label{tab:ff-data}
  \renewcommand{\arraystretch}{1.2}
  \scalebox{0.68}{
  \begin{tabular}{|l|l|l|l|l|l|l|l|l|l|l|l|l|l|}
    \cline{2-14}
      \multicolumn{1}{c|}{} & \multirow{2}{*}{\textbf{Ra}} &
      \multicolumn{3}{c|}{$v_{max}(x, 0.5)$} & \multicolumn{3}{c|}{$x$}  &
      \multicolumn{3}{c|}{$u_{max}(0.5, y)$} & \multicolumn{3}{c|}{$y$}
    \\ \cline{3-14}
      \multicolumn{1}{c|}{} & \multicolumn{1}{c|}{} & present &
      \multicolumn{1}{c|}{\cite{couturier2000performance}} &
      \multicolumn{1}{c|}{\cite{kosec2008solution}} & present &
      \multicolumn{1}{c|}{\cite{couturier2000performance}} &
      \multicolumn{1}{c|}{\cite{kosec2008solution}} & present &
      \multicolumn{1}{c|}{\cite{couturier2000performance}} &
      \multicolumn{1}{c|}{\cite{kosec2008solution}} & present &
      \multicolumn{1}{c|}{\cite{couturier2000performance}} &
      \multicolumn{1}{c|}{\cite{kosec2008solution}}
    \\ \hline \hline
      \multirow{3}{*}{\textbf{2D}} & $10^6$ & 0.2628    & 0.2604   &
      0.2627  & 0.037  & 0.038 & 0.039 & 0.0781    & 0.0765   & 0.0782  & 0.847
      & 0.851 & 0.861
    \\ \cline{2-14}
      & $10^7$ & 0.2633    & 0.2580   & 0.2579  & 0.022  & 0.023 & 0.021 & 0.0588
      & 0.0547   & 0.0561  & 0.870  & 0.888 & 0.900
    \\ \cline{2-14}
      & $10^8$ & 0.2557    & 0.2587   & 0.2487  & 0.010  & 0.011 & 0.009 & 0.0314
      & 0.0379   & 0.0331  & 0.918  & 0.943 & 0.930
    \\ \hline \hline
      \multicolumn{1}{c|}{} & \multirow{2}{*}{\textbf{Ra}} &
      \multicolumn{3}{c|}{$w_{max}(x, 0.5,0.5)$} & \multicolumn{3}{c|}{$x$}  &
      \multicolumn{3}{c|}{$u_{max}(0.5, 0.5, z)$} & \multicolumn{3}{c|}{$z$}
    \\ \cline{3-14}
      \multicolumn{1}{c|}{} & \multicolumn{1}{c|}{} & present &
      \multicolumn{1}{c|}{\cite{wang2017numerical}}
      & \multicolumn{1}{c|}{\cite{fusegi1991numerical}} & present &
      \multicolumn{1}{c|}{\cite{wang2017numerical}}
      & \multicolumn{1}{c|}{\cite{fusegi1991numerical}} & present &
      \multicolumn{1}{c|}{\cite{wang2017numerical}}
      & \multicolumn{1}{c|}{\cite{fusegi1991numerical}} & present &
      \multicolumn{1}{c|}{\cite{wang2017numerical}} &
      \multicolumn{1}{c|}{\cite{fusegi1991numerical}}
    \\  \hline
      \multirow{3}{*}{\textbf{3D}} & $10^4$ & 0.2295 & 0.2218   & 0.2252  & 0.850 & 0.887 & 0.883
      & 0.2135    & 0.1968   & 0.2013  & 0.168 & 0.179 & 0.183 \\ \cline{2-14}
      & $10^5$ & 0.2545    & 0.2442   & 0.2471  & 0.940  & 0.931 & 0.935 & 0.1564 & 0.1426 & 0.1468
      & 0.144  & 0.149 & 0.145
    \\ \cline{2-14}
      & $10^6$ & 0.2564 & 0.2556 & 0.2588 & 0.961  & 0.965 & 0.966 & 0.0841
      & 0.0816   & 0.0841 & 0.143 & 0.140 & 0.144
    \\ \hline
  \end{tabular}
  }
\end{table}

\begin{figure}[h!]
  \centering
  \includegraphics[width=0.47\linewidth]{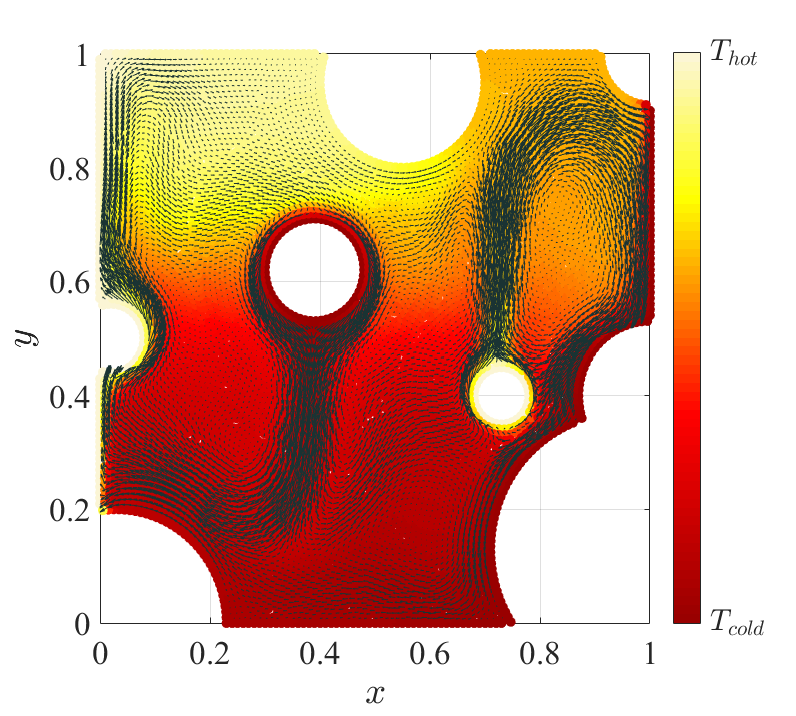}
  \raisebox{-14pt}{\includegraphics[width=0.505\linewidth]{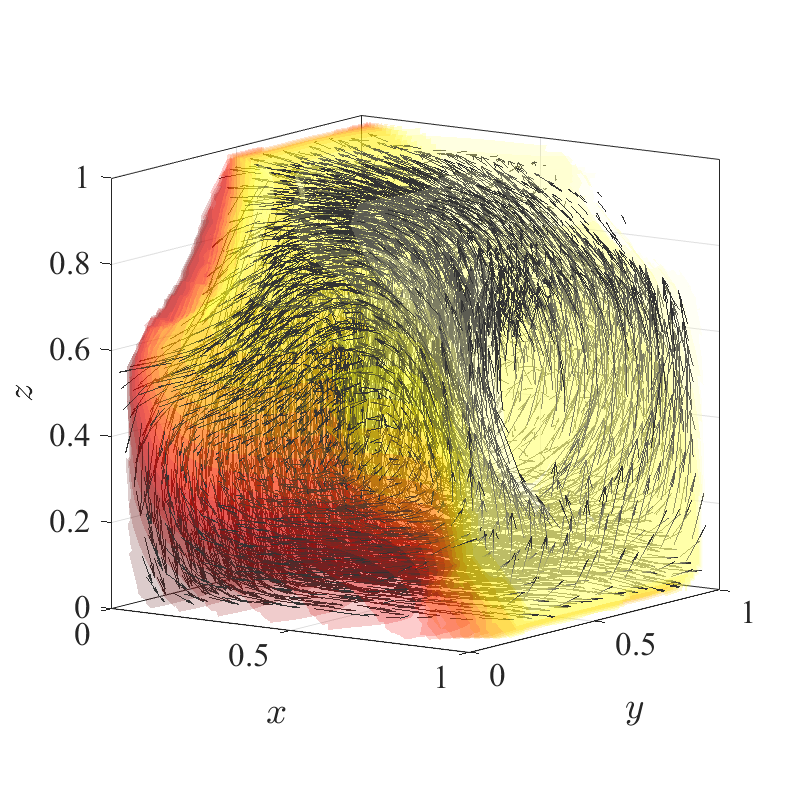}}
  \vspace{-14pt}
  \caption{Temperature contour and velocity quiver plots of solutions
   in irregular 2D domain (left) and irregular 3D domain (right).}
  \label{fig:dvd-irreg}
\end{figure}

Finally, in~\cref{fig:dvd-irreg} we demonstrate the solution of transient
convection dominated problem in an irregular 2D and 3D domain with mixed
Dirichlet-Neumann boundary conditions on nodes positioned with the proposed
algorithm. Note that this case, a solution of natural convection in an
irregular domain, includes several potential complications, such as Neumann
boundary conditions on curved boundaries, concavities, convection dominated
transport and non-linearities.

\section{Conclusions}
\label{sec:conclusions}
A new algorithm for generating variable density node distributions in interiors of
arbitrary dimensions is proposed. The algorithm has many desirable properties,
such as direction independence, support for irregular domains by only discretizing
the area actually contained in the domain interior, good compatibility with boundary
discretizations and good scaling behavior.
We prove that the time complexity of the proposed algorithm scales as
$O(N)$ for constant spacing and $O(N \log N)$ for variable spacing.
A minimal nodal spacing guarantee for constant and variable nodal spacing functions is also proven.
With examples it is shown that the proposed algorithm
produces locally smooth distributions that are suitable for RBF-FD method
for solving partial differential equations. The algorithm is compared against two
other state-of-the-art algorithms, and the summary of the findings is presented
in~\cref{tab:compare}.

\begin{table}[h]
  \centering
  \caption{Comparison of FF, SKF and PNP algorithms.}
  \label{tab:compare}
  \renewcommand{\arraystretch}{1.2}
  \scalebox{0.80}{
  \begin{tabular}{|m{4.2cm}|m{3cm}|m{3cm}|m{3.4cm}|} \hline
    \textbf{property / algorithm}   & \textbf{FF} & \textbf{SKF} & \textbf{PNP} \\ \hline  \hline
    supports variable density & yes & no & yes \\ \hline
    supports 3-D distributions & no & yes & yes \\ \hline
    supports irregular domains & yes, using BB & yes, using OBB & yes, natively \\ \hline
    compatibility with \newline boundary nodes & n/a & no & yes \\ \hline
    dimension independence & no & yes & yes \\ \hline
    direction independence & no & yes & yes  \\ \hline
    randomized & minimal \newline (only starting line) & yes (fully) & yes (controlled) \\ \hline
    minimal spacing guarantees & no & yes (constant $h$) & yes (constant \newline and variable $h$) \\ \hline
    time complexity & $O\Big(n\big(\frac{|\bb(\Omega)|}{|\Omega|}N\big)^{1.5}\Big)$
    & $O\left(n \frac{|\obb(\Omega)|}{|\Omega|} N\right)$
    & $O(nN\log N)$ \newline ($O(nN)$ if $h$ constant) \\ \hline
    computational time &
    best for smaller $N$, \newline \unit[5]{s} for $10^6$ nodes  &
    \unit[6]{s} for $10^6$ nodes &
    \unit[10]{s} for $10^6$ nodes, \newline \unit[2]{s} if $h$ constant  \\ \hline
    PDE accuracy & satisfactory & satisfactory with \newline larger support sizes & satisfactory \\ \hline
    number of free parameters & 1 (no.\ of cand.\ $n$) & 1 (no.\ of cand.\ $n$) & 1 (no.\ of cand.\ $n$) \\ \hline
  \end{tabular}
  }
\end{table}

The algorithm is also included in the Medusa library~\cite{medusa} for solving PDEs
with strong form meshless methods, but a standalone implementation of the algorithm is
available from the library's website as well~\cite{standalone}.

At least three directions are open for future research. The first one
deals with effective adaptive modification of parts of generated distributions
with target complexity $O(N_{\text{old}} + N_{\text{new}})$, where $N_\text{old}$ and
$N_{\text{new}}$ stand for the number of removed old nodes and the number of added new nodes.
The second direction is to generalize the algorithm to (parametric) surfaces,
again with desired $O(N)$ time complexity irrespective of the surface.
The third direction is to investigate parallelization opportunities on different
parallel architectures ranging from shared memory multi-core central processing units (CPUs)
and general purpose graphics processing units (\hbox{GPGPUs}) to distributed computing.
A potential approach, suitable for shared memory, is to independently build the
discretization from several seed nodes. The bottleneck in such an approach is the
manipulation of global kd-tree search structure, especially on \hbox{GPGPUS}. Alternative
simpler search structures, such as spatial grids, could be used instead,
as is common practice in computer graphics community. Second option, also suitable
for distributed computing, is via domain decomposition, where main problems arise in
load balancing and appropriate partitioning of the complex higher dimensional domains.

\section*{Acknowledgments}
The authors would like to acknowledge the financial
support of the Research Foundation Flanders (FWO), The Luxembourg National
Research Fund (FNR) and Slovenian Research Agency (ARRS) in the framework of the
FWO Lead Agency project: G018916N Multi-analysis of fretting fatigue using
physical and virtual experiments, the ARRS research core funding No.\ P2-0095
and the Young Researcher program PR-08346.

\bibliographystyle{siamplain}
\bibliography{references}
\end{document}